\documentclass[12pt]{article}
\usepackage{amsmath, amsthm}
\usepackage{amsfonts}
\usepackage{amssymb}
\usepackage[usenames,dvipsnames]{color}
\usepackage{graphicx}
\usepackage{hyperref}
\usepackage{bm}
\usepackage{listings}

\usepackage{color,calc}

\makeatletter
\definecolor{shade}{gray}{0.8}
        {
          \raggedright
        \setlength{\rightmargin}{\leftmargin}
        \setlength{\itemsep}{-12pt}
        \setlength{\parsep}{20pt}
        \begin{lrbox}{\@tempboxa}%
        \begin{minipage}{\linewidth-2\fboxsep}
        }%
        {
        \end{minipage}%
        \end{lrbox}%
        \fcolorbox{black}{shade}{\usebox{\@tempboxa}}\newline
        }%
\makeatother

\newtheorem{proposition}{Proposition}
\newtheorem{theorem}{Theorem}

\newtheorem{cor}{Corollary}

\newcommand{\iu}{\mathrm{i}} 

\newcommand{\dint}{\displaystyle\int}

\renewcommand{\eqref}[1]{\hyperref[#1]{(\ref*{#1})}}

\newcommand{\dd}{\mathrm{d}}

\newtheorem{remark}{Remark}

\newcommand*{\pref}[1]{\hyperref[#1]{(\ref*{#1})}}
\newcommand*{\refpref}[2]{\hyperref[#2]{\ref*{#1}(\ref*{#2})}}

\newcommand{\p}{\mathbb{P}}

\newcommand{\ud}{\mathrm{d}}

\newcommand{\R}{\mathbb{R}}

\definecolor{amethyst}{rgb}{0.6, 0.4, 0.8}
\definecolor{applegreen}{rgb}{0.55, 0.71, 0.0}
\definecolor{aqua}{rgb}{0.0, 1.0, 1.0}
\definecolor{asparagus}{rgb}{0.53, 0.66, 0.42}
\definecolor{amber(sae/ece)}{rgb}{1.0, 0.49, 0.0}
 	\definecolor{armygreen}{rgb}{0.29, 0.33, 0.13}
	\definecolor{shitbrown}{rgb}{0.43, 0.21, 0.1}
	\definecolor{brightpink}{rgb}{1.0, 0.0, 0.5}
	\definecolor{brightube}{rgb}{0.82, 0.62, 0.91}
	 	\definecolor{byzantine}{rgb}{0.74, 0.2, 0.64}
		\definecolor{chartreuse(web)}{rgb}{0.5, 1.0, 0.0}

\definecolor{amethyst}{rgb}{0.6, 0.4, 0.8}
\definecolor{applegreen}{rgb}{0.55, 0.71, 0.0}
\definecolor{aqua}{rgb}{0.0, 1.0, 1.0}
\definecolor{asparagus}{rgb}{0.53, 0.66, 0.42}
\definecolor{amber(sae/ece)}{rgb}{1.0, 0.49, 0.0}
 	\definecolor{armygreen}{rgb}{0.29, 0.33, 0.13}
	\definecolor{shitbrown}{rgb}{0.43, 0.21, 0.1}
	\definecolor{brightpink}{rgb}{1.0, 0.0, 0.5}
	\definecolor{brightube}{rgb}{0.82, 0.62, 0.91}
	 	\definecolor{byzantine}{rgb}{0.74, 0.2, 0.64}
		\definecolor{chartreuse(web)}{rgb}{0.5, 1.0, 0.0}

\title{Attraction to and repulsion  from   a subset of the unit sphere for isotropic stable L\'evy processes}

\author{  Andreas E. Kyprianou\thanks{Department of Mathematical Sciences, University of Bath, Claverton Down, Bath, BA2 7AY, UK. Email: \texttt{a.kyprianou@bath.ac.uk}}, \ Sandra Palau\thanks{Department of Statistics and Probability, Instituto de Investigaciones en Matem\'aticas Aplicadas y en Sistemas, Universidad Nacional Aut\'onoma de M\'exico, M\'exico.
E-mail: \texttt{sandra@sigma.iimas.unam.mx}}
\
and 
Tsogzolmaa Saizmaa\thanks{Department of Mathematical Sciences, University of Bath, Claverton Down, Bath, BA2 7AY, UK and National University of Mongolia, Baga-toiruu, Sukhbaatar district, Ulaanbaatar, Mongolia.  Email: \texttt{t.saizmaa@bath.ac.uk}}
}  
\begin{document}

\maketitle

\begin{abstract}
\noindent Taking account of recent developments in the representation of $d$-dimensional isotropic stable L\'evy  processes as self-similar Markov processes, we consider a number of new ways to condition its path. Suppose that $\mathsf{S}$ is a region of the unit sphere $\mathbb{S}^{d-1} = \{x\in \mathbb{R}^d: |x| =1\}$.   We construct the aforesaid stable L\'evy process conditioned to approach $\mathsf{S}$ continuously from either inside or outside of the sphere. Additionally, we show that 
these processes are in duality with the stable process conditioned to remain inside the sphere and absorb continuously at the origin and to remain outside of the sphere, respectively. Our results extend the recent contributions of \cite{Phil}, where  similar conditioning is considered, albeit in one dimension as well as providing analogues of the same and very classical results for Brownian motion, \cite{doob}. As in \cite{Phil}, we appeal to recent fluctuation identities related to the deep factorisation of stable processes, cf. \cite{Deep, Deep3,  KRSg}.

\medskip

\medskip

\end{abstract}

\section{Introduction}

Let $X=(X_t, t \geq 0)$ be a $d$-dimensional  stable  L\'evy process with probabilities $(\mathbb{P}_x, x \in \mathbb{R}^d)$. This means that $X$ has c\`adl\`ag paths with stationary and independent increments as well as respecting a property of  self-similarity: There is an $\alpha>0$ such that, for $c>0,$ and $x \in {\mathbb{R}^d},$
under $\mathbb{P}_x$, the law of  $(cX_{c^{-\alpha}t}, t \geq 0)$ is equal to  $\mathbb{P}_{cx}$.
It turns out that stable L\'evy processes necessarily have the scaling index $\alpha\in (0,2]$. The case $\alpha=2$ pertains to a standard $d$-dimensional Brownian motion, thus has a continuous path. The processes we construct are arguably less interesting in the diffusive setting and thus we restrict ourselves  to  the isotropic pure jump setting of $\alpha\in(0,2)$  in dimension $d\geq2$. 

\smallskip

 To be more precise, this means, for all orthogonal transformations $U:\mathbb{R}^d \mapsto \mathbb{R}^d$ and $x \in \mathbb{R}^d,$ 
\[
 \quad \textit{the law of} \quad (UX_t, t \geq 0) \textit{ under} \quad \mathbb{P}_x  \textit{ is equal to }  (X_t,t\geq 0) \textit{ under } \mathbb{P}_{Ux}.
 \]
For convenience, we will henceforth refer to $X$ just as a {\it stable process}. 
\smallskip

The stable L\'evy  process has a the jump measure $\Pi$ that satisfies 
$$\Pi(B) = \frac{2^{\alpha} \Gamma((d+\alpha)/2)}{\pi^{d/2} |\Gamma(-\alpha/2)|} \int_B \frac{1}{|y|^{\alpha+d}} {\rm d} y, \quad B \subseteq \mathcal{B}(\mathbb{R}^d).$$ 
{The constant in the definition of $\Pi(B)$ can be arbitrary, however, our choice corresponds to the one one that allows us to identify the characteristic exponent L\'evy process as}
\[
\Psi(\theta)=-\frac{1}{t}\log \mathbb{E} ({\rm e}^{\iu\theta \cdot X_t}) =|\theta|^{\alpha}, \quad \theta \in \mathbb{R}^d,
\]
where we write $\mathbb{P}$ in preference to $\mathbb{P}_0$; more precisely, the coefficient of $|
\theta|^\alpha$ is one.

\smallskip

In this article, we characterise the law of  a stable process conditioned to hit continuously a part of the surface, say $\mathsf{S}\subseteq\mathbb{S}^{d-1}= \{x\in \mathbb{R}^d: |x| =1\}$, either from the inside or from the outside of the unit sphere. We develop an expression for {the law of} the limiting point of contact on $\mathsf{S}$. Moreover, we show that, when time reversed from the strike point on $\mathsf{S}$, the resulting process can also be seen as a conditioned stable process. The extreme cases that $\mathsf{S} = \mathbb{S}^{d-1}$ (the whole unit sphere) and $\mathsf{S} = \{\vartheta\}\in \mathbb{S}^{d-1}$ (a single point on the unit sphere) are included in our analysis, however,  we will otherwise insist that the Lebesgue surface measure of $\mathsf{S}$ is strictly positive.

\smallskip

Our results relate to the recent work of \cite{Phil}, who considered a real valued L\'evy process conditioned to continuously approach the boundary of the interval $[-1,1]$ from the outside. {\it In order to avoid repetition, we always remain in two or more dimensions. } As in \cite{Phil}, we rely heavily on recent fluctuation identities that are connected to  the deep factorisation of the stable process; cf. \cite{Deep, Deep3, KRSg}. The results here are also related the classical results of Doob \cite{doob}, who deals with similar conclusions for Brownian motion and as well as echoing the general theory of conditioned stochastic processes in the potential-analytic sense (via a Doob $h$-transform), see e.g. Chapter 14 of \cite{Chung-Walsh}.

\section{Attraction towards $\mathsf{S}$}

For convenience, we will work with the definition $\mathbb{B}_{d} = \{x\in\mathbb{R}^d: |x|<1\}$. Let $\mathbb{D}(\mathbb{R}^d)$ denote the space of c\`adl\`ag paths $\omega:[0,\infty)\to\mathbb{R}^d\cup\partial$ with lifetime $\texttt{k}(\omega)=\inf\{s>0 :\omega(s)=\partial\}$, where $\partial$ is a cemetery point. The space $\mathbb{D}(\mathbb{R}^d)$ will be equipped with the Skorokhod topology, with its Borel $\sigma$-algebra $\mathcal{F}$ and natural filtration $(\mathcal{F}_t, t\geq 0)$. 
The reader will note that we will also use 
a similar notion for $\mathbb{D}(\mathbb{R}\times\mathbb{S}^{d-1})$ later on in this text in the obvious way.
We will always work with $X = (X_t,t\geq0)$ to mean  the coordinate process defined on the space $\mathbb{D}(\mathbb{R}^d)$. Hence, the notation of the introduction indicates that $\mathbb{P}= (\mathbb{P}_x, x\in\mathbb{R}^d)$ is such that $(X,\mathbb{P})$ is our stable process.

\smallskip

{Consider a subset $\mathsf{S}\subseteq\mathbb{S}^{d-1}$ such that it has strictly positive Lebesgue surface measure or it is a point.} We want to construct the law of $X$ conditioned to approach {$\mathsf{S}$}  continuously from within $ \bar{\mathbb{B}}_d^c : =\mathbb{R}^{d}\setminus \bar{\mathbb{B}}_d$. {From a potential-theoretic perspective, this law can be obtained as a Doob $h$-transform of the killed stable process in $\overline{\mathbb{B}}_d^c$, provided that $h$ is a positive harmonic function in $\overline{\mathbb{B}}_d^c$ which is equal to zero in $\overline{\mathbb{B}}_d$ and which goes to zero at infinity and at
$\mathbb{S}^{d-1}\setminus \mathsf{S}$; cf. \cite[Chapter 14]{Chung-Walsh}. In this paper, we want to give a probabilistic construction, which identifies a more physical meaning to the conditioning in terms of the paths of the stable process; see e.g. the classical work of \cite{BD, ChaumontDoney}.} Similarly, we want the law of $X$ conditioned to approach {$\mathsf{S}$} continuously from within $\mathbb{B}_d$. More precisely, via an appropriate limiting procedure, we want to build a new family of probabilities $\mathbb{P}^{\vee}= (\mathbb{P}^{\vee}_x, x\in\bar{\mathbb{B}}_d^c)$ such that 
\[
\mathbb{P}^{\vee}_x(X_s\in\bar{\mathbb{B}}_d^c, s<\texttt{k} \text{ and } X_{\texttt{k}-} \in \mathsf{S}) = 1,\qquad x\in\bar{\mathbb{B}}_d^c,
\]
with a similar statement holding when the conditioning is undertaken from within $\mathbb{B}_d$.
\smallskip

As we are considering two or higher dimensions, the process $(X,\mathbb{P})$ is transient in the sense that $\lim_{t\to\infty}|X_t| = \infty$ almost surely. Defining
\[
\underline{G}(t):=\sup\{s \le t \colon |X_s|=\inf_{u \le s} |X_u|\}, \qquad t \geq 0,
\]
we thus have by  monotonicity and the transience of $(X,\mathbb{P})$ that $\underline{G}(\infty): = \lim_{t\to\infty} \underline{G}(t)$ exists and, moreover, $X_{\underline{G}(\infty)}$ describes the point of closest reach to the origin in the range of $X$.

\smallskip

We can similarly define $\overline{G}(t)=\sup\{s \le t \colon |X_s|=\sup_{u \le s} |X_u|\}$, $t\geq 0$, so that $\overline{G}(\tau^\ominus_1-)$ is the point of 
furthest reach from the origin prior to exiting $\mathbb{B}_d$, where 
\[
\tau^\ominus_1 = \inf\{t>0: |X_t|>1\}.
\]

Let us turn to what we mean by conditioning to attract to the set $\mathsf{S}$ from either the interior or the exterior of the sphere. {If $\mathsf{S}$ is not a point,} we define 
$A_{\varepsilon}=\{r\theta : r \in (1,1+\varepsilon), \theta \in \mathsf{S}\}$ and
$B_{\varepsilon}=\{r\theta : r \in (1-\varepsilon,1), \theta \in \mathsf{S}\}$, for $0<\varepsilon<1$ and define the corresponding events  $C_\varepsilon^{\vee}:=\{X_{\underline{G}(\infty)} \in A_{\varepsilon}\},$ and  $C_\varepsilon^{\wedge}:=\{X_{\overline{G}(\tau^\ominus_1-)} \in B_{\varepsilon}\}$. 
Let 
\[
\tau^\oplus_1 = \inf\{t>0: |X_t|<1\}.
\]
We are interested in the asymptotic conditioning 
\begin{equation}
\mathbb{P}_x^\vee (A, \,t<\texttt{k})=\lim_{\varepsilon \rightarrow 0} 
\mathbb{P}_x (A, \,t<\tau^\oplus_1| C_\varepsilon^{\vee}) ,
\label{lim1}
\end{equation}
when $ x \in \bar{\mathbb{B}}_d^c $ and 
\begin{equation}
\mathbb{P}_x^\wedge (A, \,t<\texttt{k})=\lim_{\varepsilon \rightarrow 0} \mathbb{P}_x (A,\,t<\tau^\ominus_1 | C_\varepsilon^{\wedge})  ,
\label{lim2}
\end{equation}
when $x \in \mathbb{B}_d$, 
for all $A\in \mathcal{F}_t$.

\smallskip

When $\mathsf{S} = \{\vartheta\}\in \mathbb{S}^{d-1}$, we need to adapt slightly the sets $A_\varepsilon$ and $B_\varepsilon$ so that 
$A_\varepsilon = \{r\phi \colon r \in (1,1+\varepsilon), \phi\in\mathbb{S}^{d-1}, |\phi-\vartheta|<\varepsilon\}$ and $B_\varepsilon = \{r\phi\colon r \in (1-\varepsilon,1),  \phi\in\mathbb{S}^{d-1}, |\phi-\vartheta|<\varepsilon\}.$

\smallskip

We will go a little further in due course and give a fuller description of these two conditioned processes  by including the cases that  $X$ is issued from the unit sphere itself but not within $\mathsf{S}$, i.e. $\mathbb{S}^{d-1}\setminus\mathsf{S}$. For now, we have our first main result, given immediately below, for which we define the function

\begin{equation}
H_\mathsf{S}(x)= \left\{
\begin{array}{ll}
 ||x|^2-1|^{\alpha/2}\dint_\mathsf{S} |\theta-x|^{-d} \sigma_1({\rm d}\theta)& \text{ if }\sigma_1(\mathsf{S})>0,\\
&\\
||x|^2-1|^{\alpha/2}|\vartheta-x|^{-d}  &\text{ if }\mathsf{S} =\{\vartheta\},
\end{array}
\right.
\label{H_S}
\end{equation}
for $|x|\neq 1$, where $\sigma_1({\rm d}\theta)$ is the Lebesgue surface measure on $\mathbb{S}^{d-1}$ normalised to have unit mass. 
It is worthy of note that, when $\mathsf{S} = \mathbb{S}^{d-1}$, the integral in \eqref{H_S}
 can be computed precisely. 
Indeed, up to an unimportant (for our purposes) multiplicative constant, $C>0$, which may change from line to line, we note that, for $|x|>1$,
\begin{align*}
\int_{\mathbb{S}^{d-1}}|x- \theta|^{-d}\sigma_1(\dd \theta)& = C\int_0^\pi\frac{(\sin\phi)^{d-2}}{(|x|^2 -2 |x|\cos\phi +1 )^{d/2}}\dd \phi\\
&= C|x|^{-d}{_2}F_1(d/2, 1; d/2, |x|^{-2})\\
&=C |x|^{-d} \left(1 - \frac{1}{|x|^2}\right)^{-1},
\end{align*}
where we have used the hypergeometric identity in \eqref{3.665} of the Appendix. We can perform a similar calculation when $|x|<1$ and, obtain, up to a multiplicative constant, $C>0$, that 
\begin{equation}
\int_{\mathbb{S}^{d-1}}|x- \theta|^{-d}\sigma_1(\dd \theta)= C \left(1 - {|x|^2}\right)^{-1} .
\label{lessthanone}
\end{equation}
All together, noting that we may ignore multiplicative constants, we have 
\begin{equation}
H_{\mathbb{S}^{d-1}}(x)= \left\{
\begin{array}{ll}
|x|^{\alpha-d}\left(1 - |x|^{-2}\right)^{\frac{\alpha}{2}-1}& \text{ if }|x|>1\\
&\\
(1-|x|^2)^{\frac{\alpha}{2}-1} &\text{ if }|x|<1.
\end{array}
\right.
\label{justH}
\end{equation}

\smallskip

{As the next result will make clear, $H_\mathsf{S}$ is a positive harmonic function for both $(X_t, t<\tau^\ominus_1)$ and $(X_t, t<\tau^\oplus_1)$.  From the potential-theoretic perspective, it can be described as an integral of the Martin kernel over $\mathsf{S}$. Then, by the  Martin boundary theory, the $h$-conditioned process will approach $\mathsf{S}$ with probability one, see \cite[Chapter 14]{Chung-Walsh} as well as the classical results of Doob  for Brownian motion, cf. Theorem 7.1 \cite{doob}.}

\begin{theorem}[Stable process conditioned to attract to $\mathsf{S}$ continuously from one side]\label{main1}
Let $\mathsf{S} \subseteq\mathbb{S}^{d-1}$ be an {closed} set with $\sigma_1(\mathsf{S})>0$ or $\mathsf{S}=\{\vartheta\}$ for a fixed point $\vartheta \in \mathbb{S}^{d-1}.$ Then for all points of issue $x \in \mathbb{R}^d \setminus \mathbb{S}^{d-1}$ we have
\begin{equation}
\left.\frac{d \mathbb{P}_x^\vee}{d \mathbb{P}_x} \right| _{\mathcal{F}_t} = \mathbf{1}_{(t < {\tau^\oplus_1})} \frac{H_\mathsf{S}(X_t)}{H_\mathsf{S}(x)},\qquad \text{ if }x\in\bar{\mathbb{B}}_d^c
\label{out}
\end{equation} 
and otherwise 
\begin{equation}
\left.\frac{d \mathbb{P}_x^\wedge}{d \mathbb{P}_x} \right| _{\mathcal{F}_t} = \mathbf{1}_{(t < {\tau^\ominus_1})} \frac{H_\mathsf{S}(X_t)}{H_\mathsf{S}(x)},\qquad \text{ if }x\in \mathbb{B}_d.
\label{in}
\end{equation} 
In particular, $(\mathbb{P}_x^\vee, x\in \bar{\mathbb{B}}_d^c)$ and $(\mathbb{P}_x^\wedge, x\in \bar{\mathbb{B}}_d^c) $ are Markovian families.
\end{theorem}
\begin{remark}\rm
The choice of limiting conditioning procedure that we have used reflects a similar approach taken in \cite{Phil} in one dimension.
It is worth noting at this point that the choice of $C_\varepsilon^{\vee}$ and $C_\varepsilon^{\wedge}$ are by no means the only possibilities as far as  performing a limiting   conditioning that results in \eqref{out} and \eqref{in}. For example, once the reader is familiar with the proof of Theorem \ref{main1}, it will quickly become clear that, when $\mathsf{S}$ is not a singleton, by defining e.g.  
$C_\varepsilon^{\vee}= \{X_{\tau^\oplus_1}  \in B_\varepsilon\}$, or indeed $C_\varepsilon^{\vee}= \{X_{\tau^\oplus_1 -}  \in A_\varepsilon\}$, 
 the limit \eqref{lim1}
will still produce the change of measure \eqref{out}. Once the reader is familiar with the proof of Theorem \ref{main1}, it is a worthwhile exercise to verify the two proposed alternative definitions of $C_\varepsilon^{\vee}$ for the limiting process  by appealing to the fluctuation identities in e.g. \cite{Deep3}. Other definitions of $C_\varepsilon^{\vee}$ giving a consistent limit may indeed also be possible. 
\end{remark}

Whilst the above theorem deals with the construction of the conditioned process up to but not including its terminal position, we characterise the latter in the next result, which resonates with Theorem 14.8 of \cite{Chung-Walsh}.  
\begin{proposition}[Distribution of the hitting location]\label{prop2}Suppose that $\mathsf{S} \subseteq\mathbb{S}^{d-1}$ be a {closed} set with $\sigma_1(\mathsf{S})>0$.
Let $\mathsf{S}' $ be an {closed} subset of $\mathsf{S}$. Then for any $x \in \mathbb{R}^d \setminus \bar{\mathbb{B}}_d,$ we have 
\begin{equation}
\mathbb{P}_x^{\vee} (X_{\emph{\texttt{k}}-} \in\mathsf{S}' ) = \frac{\int_{\mathsf{S}'} |\theta-x|^{-d}\sigma_1({\rm d}\theta) }{\int_\mathsf{S} |\theta-x|^{-d}\sigma_1({\rm d}\theta) },
\label{hitdist}
\end{equation}
with an identical result holding for $X_{\emph{\texttt{k}}-}$ under $\mathbb{P}^{\wedge}_x$, with $x\in\mathbb{B}_d$.
\end{proposition}

\section{Lamperti--Kiu representation and radial excursions}\label{radialexcursionsection}
The basic definition of the stable process conditioned to attract continuously to $\mathsf{S}$ from one side is not quite complete. Strictly speaking, we could think about defining the process to include the points of issue in $\mathbb{S}^{d-1}\setminus \mathsf{S}$. It turns out that this is possible. However, we first need to remind the reader of the recently described radial excursion theory, see \cite{Deep3, cone}. The starting point for the aforementioned is the Lamperti--Kiu transform which identifies the stable process as a self-similar Markov process. 

\smallskip
To describe it, we need to introduce the notion of a Markov Additive Process, henceforth written MAP for short.
Let $\mathbb{S}^{d-1} =\{x\in\mathbb{R}^d: |x|=1\}$. With an abuse  of previous notation, we say that $(\Xi,\Upsilon) =( (\Xi_t, \Upsilon_t), t\geq 0)$ is a MAP if it is Feller process on $\mathbb{R}{^n}\times\mathbb{S}^{d-1},$ with probabilities $\texttt{P}_{x,\theta}$, $x\in\mathbb{R}{^n}$, $\theta\in\mathbb{S}^{d-1}$,  such that, for any $t\geq 0$, the conditional law of the process $((\Xi_{s+t}-\Xi_t,\Upsilon_{s+t}):s\geq 0)$, given  $((\Xi_u,\Upsilon_u), u\leq t),$ is that of $(\Xi,\Upsilon)$ under $\texttt{P}_{0,\theta}$, with $\theta=\Upsilon_t$. For a MAP pair $((\Xi_t,\Upsilon_t), t\geq 0)$, we call $\Xi$ the {\it ordinate} and $\Upsilon$ the {\it modulator}.

\smallskip

According to one of the main results in \cite{ACGZ}, there exists a MAP on ${\mathbb{R}\times\mathbb{S}^{d-1}}$,  which we will henceforth write as $(\xi, \Theta)$, with probabilities $\mathbf{P} =(\mathbf{P}_{x,\theta}, x\in\mathbb{R}, \theta\in\mathbb{S}^{d-1})$ such that the $d$-dimensional  stable process can be written
   \begin{equation}\label{eq:lamperti_kiu}
    X_t  = \exp\{\xi_{\varphi(t)}\} \Theta_{\varphi(t)}\qquad t \geq 0,
  \end{equation}
  where
 \begin{equation}
 \label{varphi}
 \varphi(t) = \inf\left\{s>0: \int_0^s{\rm e}^{\alpha\xi_u}\dd u>t\right\}.
 \end{equation}
 
Whilst $\Theta$ alone is a Feller process, it is not necessarily true that $\xi$ alone is. However, it is a consequence of isotropy that this is the case here. Moreover, $\xi$ alone is a L\'evy process whose characteristic exponent is known (but not important in the current context); see for example \cite{CPP}. What is important for our purposes is to note for now that it has paths of unbounded variation, and therefore is regular for the upper and lower half line (in the sense of Definition 6.4 of \cite{kbook}).

 \smallskip
 
 It is not difficult to show that the pair $((\xi_t- \underline{\xi}_t, \Theta_t), t\geq 0)$, forms a strong Markov process, where $\underline\xi_t: = \inf_{s\leq t}\xi_s$, $t\geq 0$ is the running minimum
 of $\xi$.   On account of the fact that $\xi$, alone, is a L\'evy process, $(\xi_t- \underline{\xi}_t,\ t\geq0)$ is also a strong Markov process. Suppose we denote by $\ell = (\ell_t, t\geq 0)$ the local time at zero of $\xi - \underline\xi$, then we can introduce the following processes
  	\[
  	H^-_t =-\xi_{\ell^{-1}_t}\text{ and } \Theta^-_t = \Theta_{\ell^{-1}_t},\qquad t\geq 0,
  	\]
	and define $(H^-_{\ell^{-1}_t}, \Theta^-_{\ell^{-1}_t}) = (\partial,\dagger)$, a cemetery state, if $\ell^{-1}_t = \infty$.
  	Then, the pair $(\ell^{-1}, H^-)$, without reference to the associated modulation $\Theta^-$, are Markovian and play the role of the descending ladder time and height subordinators of $\xi$. Moreover, the strong Markov property tells us that $(\ell^{-1}_t, H^-_t, \Theta^-_t)$, $t\geq 0$, defines a Markov Additive Process { on $\mathbb{R}^2\times\mathbb{S}^{d-1}$}, whose first two elements are ordinates that are non-decreasing. In this sense, $\ell$ also serves as an adequate choice for the local time of the Markov process $(\xi - \underline\xi, \Theta)$ on the set $\{0\}\times\mathbb{S}^{d-1}$. 
 {(See \cite{Deep3}).}

 \smallskip
 
Suppose we define ${\texttt g}_t =\sup\{s< t: \xi_s = \underline\xi_s\},$ and recall that the regularity of $\xi$ for $(-\infty, 0)$ and $(0,\infty)$ ensures that it is well defined, as is ${\texttt g}_\infty = \lim_{t\to\infty}{\texttt g}_t$. Set 
\[
{\texttt d}_t = \inf\{s> t: \xi_s = \underline\xi_s\}.
\]
For all $t>0$ such that ${\texttt d}_t > {\texttt g}_t$ the process 
\[
(\epsilon_{{\texttt g}_t}(s), \Theta^\epsilon_{{\texttt g}_t}(s)): = (\xi_{{\texttt g}_t +s}-\xi_{{\texttt g}_t}, \Theta_{{\texttt{g}_t + s}}), \qquad s\leq \zeta_{{\texttt g}_t}: = {\texttt d}_t-{\texttt g}_t,
\]
codes the excursion of $(\xi-\underline\xi, \Theta)$ from the set $(0,\mathbb{S}^{d-1})$ which straddles time $t$. Such excursions live in the space $\underline{\mathbb{U}}(\mathbb{R}\times\mathbb{S}^{d-1})$, the space of c\`adl\`ag paths 
in $\mathbb{R}\times\mathbb{S}^{d-1}$, written in canonical form
\[
(\epsilon, \Theta^\epsilon) = ((\epsilon(t), \Theta^\epsilon(t)): t\leq \zeta)\text{ 
with lifetime  }\zeta = \inf\{s>0: \epsilon(s) <0\},
\] 
such that $(\epsilon(0),\Theta^\epsilon(0))\in \{0\}\times\mathbb{S}^{d-1}$, $(\epsilon(s), \Theta^\epsilon(s))\in (0,\infty)\times\mathbb{S}^{d-1}$, for $0<s<\zeta$,  and $\epsilon(\zeta)\in (-\infty,0]$.
\smallskip

Taking account of the Lamperti--Kiu transform \eqref{eq:lamperti_kiu}, it is natural to consider how the excursion of $(\xi - \underline{\xi}, \Theta)$ from $\{0\}\times\mathbb{S}^{d-1}$ translates into a radial excursion theory for the process 
\[
Y_t  : = {\rm e}^{\xi_t}\Theta_t, \qquad t\geq 0.
\]
Ignoring the time change in \eqref{eq:lamperti_kiu}, we see that the radial minima of the process $Y$ agree with the radial minima of the stable process $X$.
Indeed, each excursion of $(\xi - \underline{\xi}, \Theta)$ from $\{0\}\times\mathbb{S}^{d-1}$  is uniquely associated to exactly one excursion of $(Y_t/\inf_{s\leq t}|Y_s|, t\geq 0)$, from $\mathbb{S}^{d-1}$, or equivalently an excursion of $Y$ from its running radial infimum. Moreover, we see that, for  all $t>0$ such that ${\texttt d}_t > {\texttt g}_t$, 
\[
Y_{\texttt{g}_t + s} = {\rm e}^{\xi_{\texttt{g}_t}} {\rm e}^{\epsilon_{\texttt{g}_t}(s)} \Theta^\epsilon_{{\texttt g}_t}(s) = |Y_{\texttt{g}_t}|{\rm e}^{\epsilon_{\texttt{g}_t}(s)} \Theta^\epsilon_{{\texttt g}_t}(s)
, \qquad s\leq \zeta_{{\texttt g}_t}.
\]
This will be useful to keep in mind for  the forthcoming excursion computations.

\smallskip

For $t>0,$ let $R_t=\texttt{d}_t-t,$ and define the set $\mathbb{G}= \{t>0: R_{t-}=0, R_{t}>0\} = \{{\texttt g}_s: s\geq 0\}$. 
The classical theory of exit systems in  \cite{Maison} (see {Theorems (4.1) and (6.3)} therein) now implies that there exists 
 an additive functional $(\Lambda_t, t\geq 0)$   and a family of {\it excursion measures}, $(\underline{\mathbb{N}}_{\theta }, \theta\in\mathbb{S}^{d-1})$
such that:  
\begin{itemize}
\item[(i)] $\Lambda$ is  an additive functional of $(\xi,\Theta)$, has a bounded $1$-potential  and is carried by the set of times $\{t\geq0: (\xi_t - \underline{\xi}_t, \Theta_t)\in\{0\}\times\mathbb{S}^{d-1}\}$,
\item[(ii)] the map $\theta\mapsto \underline{\mathbb{N}}_{\theta}$ is an $\mathbb{S}^{d-1}$-indexed kernel on $\underline{\mathbb{U}}(\mathbb{R}\times\mathbb{S}^{d-1})$  such that $\underline{\mathbb{N}}_{\theta}(1-{\rm e}^{-\zeta})<\infty$;
\item[(iii)]we have the {\it exit formula} \begin{align}
&\mathbf{E}_{x,\theta}\left[\sum_{g\in\mathbb{G} }F((\xi_s, \Theta_s): s<g)H((\epsilon_{g}, \Theta^\epsilon_g))\right]\notag\\
&\hspace{2cm}=\mathbf{E}_{x,\theta}\left[\int_0^\infty F((\xi_s, \Theta_s): s< t)\underline{\mathbb{N}}_{ \Theta_t}(H(\epsilon, \Theta^\epsilon)){\rm d}\Lambda_t\right],
\label{exitsystem1}
\end{align}
for $x\neq 0$, where $F$ is continuous on the space of c\`adl\`ag paths $\mathbb{D}(\mathbb{R}\times\mathbb{S}^{d-1})$ and $H$ is measurable on the space of c\`adl\`ag paths $\underline{\mathbb{U}}(\mathbb{R}\times\mathbb{S}^{d-1});$
\item[(iv)] under any measure $\underline{\mathbb{N}}_{\theta}$ the process $((\epsilon(s), \Theta^\epsilon(s)), s<\zeta)$ is  a strong Markov process with the same semigroup as $(\xi, \Theta)$ killed at its first hitting time of $(-\infty,0]\times\mathbb{S}^{d-1}.$   
\end{itemize}
The couple $(\Lambda, \underline{\mathbb{N}}_{\cdot})$ is called an exit system. Note that in Maisonneuve's original formulation, the pair $\Lambda$, $\underline{\mathbb{N}}_\cdot: = (\underline{\mathbb{N}}_\theta, \theta\in\mathbb{S}^{d-1} )$ is not unique, but once $\Lambda$ is chosen the measures $(\underline{\mathbb{N}}_\theta, \theta\in \mathbb{S}^{d-1})$ exist however, are only unique up to $\Lambda$-{\it neglectable} sets, i.e. sets $\mathcal{A}$ 
such that $\mathbf{E}_{x,\theta}(\int_{0}^\infty\mathbf{1}_{\{\Theta_s\in\mathcal{A}\}}{\rm d}\Lambda_s)=0$. Another example of where this theory has been used is in the construction of excursions from a set is that  of Brownian motion away from a hyperplane; see \cite{Burdzy}.

\smallskip

Now referring back to $\ell$, the local time of $\xi -\underline\xi$ at 0, since it is an additive functional with a bounded $1$-potential, there is an exit system which corresponds to  $(\ell, \underline{\mathbb{N}}_{\cdot})$. With this choice of $\ell$ we assume that the choice of 
 $\underline{\mathbb N}_\cdot$ is fixed despite the fact that we can induce subtle variations in  $\underline{\mathbb N}_\cdot$ on  a $\Lambda$-negligible set of $\theta\in\mathbb{S}^{d-1}$ e.g. by setting $\underline{\mathbb N}_\theta \equiv 0$ there. The reader is referred to Chapter VII of \cite{Blum} for further discussion on this matter.
%
Note that  $\underline{\mathbb N}_\theta$ is not isotropic in $\theta$. For example, excursions that begin at the `North Pole', say $\mathtt{1}$, are, with high frequency, arbitrarily small and hence will end near to $\texttt{1}$.  That said, depending on the event $A$, it is possible that $\underline{\mathbb N}_\theta(A)$ does not depend on $\theta\in\mathbb{S}^{d-1}$; for example,  $\underline{\mathbb{N}}_{\theta}(\zeta=\infty)$. The reason for this is that it must agree with the rate at which the infinite excursion of $\xi-\underline{\xi}$ occurs, according to the local time $\ell$. More generally, we have that, for all orthogonal transformations $U:\mathbb{R}^d \mapsto \mathbb{R}^d$ and $f$ such that $\underline{\mathbb N}_\theta(f(\epsilon, \Theta^\epsilon))<\infty$, $\theta\in\mathbb{S}^{d-1}$, isotropy implies that $\underline{\mathbb N}_\theta(f(\epsilon, U\Theta^\epsilon)) = \underline{\mathbb N}_{U\theta}(f(\epsilon, \Theta^\epsilon))$.
On account of the fact that $\ell$ is only defined up to a multiplicative constant, we can use the common value of  $\underline{\mathbb{N}}_{\theta}(\zeta=\infty)$ to fix a normalisation the local time, or equivalently, of the excursion measures $(\underline{\mathbb N}_\theta, \theta\in\mathbb{S}^{d-1})$. We thus fix it to take the value of unity. The place at which this choice of normalisation becomes relevant is when we cite certain identities from (cf. \eqref{upNint} below) from \cite{Deep3}, in which this assumption was also made.  
 Henceforth, this is the exit system we will work with and the system of excursion associated to it is what we call our {\it radial excursion theory}.

\smallskip

Later in our proofs we will use a variant of the above excursion theory based on the MAP $(\overline\xi - \xi, \Theta)$, where $\overline\xi$ is the process $\overline\xi_t = \sup_{s\leq t}\xi_s$, $t\geq 0$. We leave the details until that point in the text.
With our excursion theory in hand, we can now proceed to identify the completion of Theorem \ref{main1}. 
\begin{theorem}\label{main3} Let $\mathsf{S} \subseteq\mathbb{S}^{d-1}$ be an {closed} set with $\sigma_1(\mathsf{S})>0$ or $\mathsf{S}=\{\vartheta\}$ for a fixed point $\vartheta \in \mathbb{S}^{d-1}.$  The processes $(X,\p^\vee)$ and $(X,\p^\wedge)$ can be extended in a consistent way  to include points of issue $x\in \mathbb{S}^{d-1}\setminus{\mathsf{S}}$ with pathwise continuous entry via 

\begin{equation}
    \mathbb{P}_x^\vee (X_t \in {\rm d} y, \, t<\emph{\texttt{k}}) := 
   \frac{\Gamma(d/2) }{\Gamma(\alpha/2+1)\Gamma((d-\alpha)/2)} \frac{H_\mathsf{S}(y)}{h(x)} \underline{\mathbb{N}}_x \left(X^\epsilon(t)\in {\rm d}y, \, t<\varsigma \right), \qquad |y| >1
    \label{veeS}
\end{equation}
and
\begin{equation}
  \mathbb{P}_x^\wedge (X_t \in {\rm d} y, \, t<\emph{\texttt{k}}): =
  \frac{\Gamma(d/2) }{\Gamma(\alpha/2+1)\Gamma((d-\alpha)/2)} \frac{H_\mathsf{S}(y)}{h(x)}    \overline{\mathbb{N}}_x \left(X^\epsilon(t) \in {\rm d} y,\, t<\varsigma \right), \qquad |y|<1,
\end{equation}
where, 
\[
h(x) = \int_\mathsf{S}|x -\theta|^{-d}\sigma_1(\dd\theta)
\]
and, for $(\epsilon,\Theta^\epsilon)$ selected from $\underline{\mathbb{U}}(\mathbb{R}\times\mathbb{S}^{d-1})$ or $\overline{\mathbb{U}}(\mathbb{R}\times\mathbb{S}^{d-1})$, respectively, 
\begin{equation}
X^\epsilon(t) = {\rm e}^{\epsilon({\varphi(t)})}\Theta^\epsilon({\varphi(t)})\text{ and }\varsigma = \varphi^{-1}(\zeta) = \int_0^\zeta |X^\epsilon(u)|^\alpha\dd u.
\label{Xexcursions}
\end{equation}
Here, pathwise continuous entry means that 
\begin{equation}
\mathbb{P}_x^\vee(\lim_{t\to0} X_t = x)= \mathbb{P}_x^\wedge(\lim_{t\to0} X_t = x)=1
\label{forcesnormalisation}
\end{equation}
 for all $x\in\mathbb{S}^{d-1}\setminus \mathsf{S}$. 
\end{theorem}
{Note, referring to the discussion preceding Theorem \ref{main3} that pertains to the choice of excursion measures and local time, given the choice of local time $\ell$ leaves a free choice of multiplicative constant in the definition of local time, which may depend on $x\in \mathbb{S}^{d-1}\setminus\mathsf{S}$. In the proof of Theorem \ref{main3}, we use a method of continuity of resolvents to pin down the aforesaid constants. We also note that extending the notion of a Doob $h$-transformed process to include certain `boundary points' in the way we have seen in Theorem \ref{main3} can be seen in e.g. \cite{Panti, ChaumontDoney} as well as the classical work of Doob \cite{doob}.}

\section{Repulsion and duality}
In this section, we want to introduce two new processes, which will turn out to be dual to $(X,\mathbb{P}^\vee)$ and $(X,\mathbb{P}^\wedge)$ in the sense of time reversal. The two processes we are interested give meaning to the stable process conditioned to remain in $\bar{\mathbb{B}}^c_d$ and  $\mathbb{B}_d$, respectively, in an appropriate sense. 
\smallskip

An important tool that we will make use of in analysing the aforesaid time reversed processes comes through the so-called 
Riesz--Bogdan--\.Zak transform, which relates path behaviour of the stable process outside of the unit sphere to its behaviour inside the unit sphere. In order to state it, we need to introduce the process $(X,\mathbb{P}^\circ)$, where the probabilities $\mathbb{P}^\circ = (\mathbb{P}^\circ_x, x\neq 0)$ are given by 
\begin{equation}
\left.\frac{\dd \mathbb{P}^\circ_x}{\dd \mathbb{P}_x}\right|_{\mathcal{F}_t} = \frac{|X_t|^{\alpha-d}}{|x|^{\alpha-d}}, \qquad \text{ on }t<\tau_\varepsilon : =\inf\{t>0 : |X_t|<\varepsilon\}
\label{COM}
\end{equation}
for all $\varepsilon>0$.
Since $\alpha<2\leq d$, we note that the change of measure rewards paths that approach the origin and punishes paths that wander far from the origin. Intuitively, it is clear that $(X,\mathbb{P}^\circ)$  describes the stable process conditioned to continuously approach the origin. Nonetheless, this heuristic can be made into a rigorous statement, see for example \cite{KALEA, Deep3, cone, KRSg}. The reader will also note from these references (and it is easy to prove that) that $(X,\mathbb{P}^\circ)$ is also a self-similar Markov process with the same index of self-similarity as $(X,\mathbb{P})$.

\begin{theorem}[Riesz--Bogdan--\.Zak transform]\label{RBSthrm}
Suppose we write $Kx = x/|x|^2$, $x\in\mathbb{R}^d$ for the classical inversion of space through the sphere $\mathbb{S}^{d-1}$. Then, in dimension $d\geq 2$,   for $x\neq 0$, $(KX_{\eta(t)}, t\geq 0)$ under $
\mathbb{P}_{x}$ is equal in law to   $(X,\mathbb{P}^\circ_{Kx})$, where $\eta(t) = \inf\{s>0: \int_0^s|X_u|^{-2\alpha}\ud u>t\}$. 
\end{theorem}

Let us return to our duality concerns.
To this end, let us introduce the probabilities 
\begin{equation}
    H^\ominus(x)=\mathbb{P}_x (\tau_1^\oplus = \infty) = \frac{\Gamma(d/2)}{\Gamma((d-\alpha)/2) \Gamma(\alpha/2)} \int_0^{|x|^2-1} (u+1)^{-d/2} u^{\alpha/2-1} \dd u,
    \label{z=oo}
\end{equation}
for $|x|>1,$ where the second inequality is lifted from \cite{Blu}, and,
\[
H^\oplus (x)= |x|^{\alpha -d} H^\ominus(Kx),
\]
for $|x|<1$.

These two functions are positive harmonic   for $X$ and 
can be used to define the two families of probabilities $\mathbb{P}^\ominus = (\mathbb{P}^\ominus_x, |x|>1)$ and $\mathbb{P}^\oplus = (\mathbb{P}^\oplus_x, |x|<1)$ via the {Doob $h$-transforms},
\begin{equation}
\left.\frac{\dd \mathbb{P}^\ominus_x}{\dd \mathbb{P}_x}\right|_{\mathcal{F}_t} = \frac{H^\ominus(X_t)}{H^\ominus(x)}\mathbf{1}_{(t<\tau^\oplus)}, 
\qquad t\geq 0, |x|>1
\label{upCOM}
\end{equation}
and, 
\begin{equation}
\left.\frac{\dd \mathbb{P}^\oplus_x}{\dd \mathbb{P}_x}\right|_{\mathcal{F}_t} = \frac{H^\oplus(X_t)}{H^\oplus(x)}\mathbf{1}_{(t<\tau^\ominus)},\qquad t\geq 0, |x|<1.
\label{downCOM}
\end{equation}
The first of these two changes of measure corresponds to  the stable process conditioned to avoid entering $\mathbb{B}_d$ by a simple restriction on the probability space (remembering that $\lim_{t\to\infty}|X_t| = \infty$). {Note} from Theorem \ref{RBSthrm} that 
\[
    H^\ominus(Kx)=\mathbb{P}_{Kx} (\tau_1^\oplus = \infty)=\mathbb{P}^\circ_{x} (\tau^{\{0\}} <\tau_1^\ominus),
\]
where $\tau^{\{0\}} = \inf\{t>0: |{X_{t-}} |= 0\}$.
The second change of measure, \eqref{downCOM}, is a composition of conditioning the stable process to be {absorbed} continuously at the origin, followed by conditioning it not to exit $\mathbb{B}_d$ via a simple restriction on the probability space (noting that  $\lim_{t\to\infty}|X_t| = 0$ under $\mathbb{P}^\circ$).

\smallskip

The reader will also note that the Riesz-Bogdan-\.Zak transform also implies a similar spatial inversion and time change must hold for the pair $(X, \mathbb{P}^\ominus)$ and $(X,\mathbb{P}^\oplus)$. 
\begin{cor}
 For $|x|>1$, $(KX_{\eta(t)}, t\geq 0)$ under $
\mathbb{P}^\ominus_{x}$ is equal in law to   $(X,\mathbb{P}^\oplus_{Kx})$, where $\eta(t) = \inf\{s>0: \int_0^s|X_u|^{-2\alpha}\ud u>t\}$.  Similarly, for $|x|<1$,  $(KX_{\eta(t)}, t\geq 0)$ under $
\mathbb{P}^\oplus_{x}$ is equal in law to   $(X,\mathbb{P}^\ominus_{Kx})$.
\end{cor}
\noindent{\bf Proof.} Suppose that $F(X_s,s\leq t)$ is a bounded $\mathcal{F}_t$-measurable function for each $t\geq 0$.  Then,  for $|x|>1$, appealing to Theorem \ref{RBSthrm}, we have 
\begin{align*}
\mathbb{E}^\ominus_x\left[F(KX_{\eta(s)}, s\leq t)\right]&= \mathbb{E}_x\left[F(KX_{\eta(s)}, s\leq t)\frac{H^\ominus(K(KX_{\eta(t)}))}{H^\ominus(x)}\mathbf{1}_{(\eta(t)<\tau^\oplus)}\right]\\
&= \mathbb{E}^\circ_{Kx}\left[F(X_{s}, s\leq t)\frac{H^\ominus(KX_{t})}{H^\ominus(x)}\mathbf{1}_{(t<\tau^\ominus)}\right]\\
&= \mathbb{E}_{Kx}\left[F(X_{s}, s\leq t)\frac{|X_t|^{\alpha -d}}{|Kx|^{\alpha -d}}\frac{H^\ominus(KX_{t})}{H^\ominus(K(Kx))}\mathbf{1}_{(t<\tau^\ominus)}\right]\\
&= \mathbb{E}_{Kx}^\oplus\left[F(X_{s}, s\leq t)\right].
\end{align*}
This shows the first half of the claim. The second part of the claim is proved using the same technique and the details are omitted for brevity given how straightforward they are.
\hfill$\square$
\smallskip

In the spirit of other cases of conditionings from an extreme boundary point (e.g. conditioning a L\'evy process  to avoid the origin, cf. \cite{Panti}, or to stay positive, cf. \cite{ChaumontDoney}), we can extend the definitions given in \eqref{upCOM} and \eqref{downCOM} by appealing to the Markov property of the excursion measures $\underline{\mathbb N}_x$ and $\overline{\mathbb N}_x$, $x\in\mathbb{S}^{d-1}$.
\begin{theorem}
\label{mainx} The processes $(X,\mathbb{P}^\ominus)$ and $(X,\mathbb{P}^\oplus)$ can be extended in a consistent way to include points of issue on $\mathbb{S}^{d-1}$. Specifically, 
\begin{equation}
\mathbb{P}^\ominus_x(X_t \in {\rm d}y) = H^\ominus(y)\underline{\mathbb N}_x\left(X^\epsilon(t) \in {\rm d}y, \, t<\varsigma \right), \qquad {x\in \mathbb{S}^{d-1} }, |y|>1
\label{Nominus}
\end{equation}
 and similarly 
 \begin{equation}
\mathbb{P}^\oplus_x(X_t \in {\rm d}y) = H^\oplus(y)\overline{\mathbb N}_x\left(X^\epsilon(t) \in {\rm d}y, \, t<\varsigma \right), \qquad {x\in \mathbb{S}^{d-1}}, |y|<1,
\label{Noplus}
\end{equation}
(specifically, the normalisation of the excursion measure is unity in both cases)
where we have used the notation given in \eqref{Xexcursions}. {As in Theorem \ref{main3},  there is pathwise continuous  entry.}
\end{theorem}
\smallskip

Our objective is to pair up  $(X, \mathbb{P}^\vee)$, $(X, \mathbb{P}^\ominus)$ and $(X, \mathbb{P}^\wedge)$, $(X, \mathbb{P}^\oplus)$ via Nagasawa's duality theorem for time reversal; cf \cite{Naga}. To this end we need to introduce the notion of $L$-times.
\smallskip

 Suppose that  $Y = (Y_t, t\leq \zeta)$ with probabilities ${\rm\texttt{P}}_x$, $x\in E$, is a regular Markov process on an open domain $E \subseteq \mathbb{R}^d$ (or more generally, a locally compact Hausdorff space with countable base), with cemetery state $\Delta$ and killing time $\zeta=\inf\{t>0: Y_t = \Delta\}$. Let us additionally write   ${\texttt P}_\nu = \int_{E}\nu(\dd a){\texttt P}_a$, for any probability measure $\nu$ on the state space of $Y$.
	
	\smallskip
	
	Suppose that $\mathcal{G}$ is the $\sigma$-algebra generated by $Y$ and  write $\mathcal{G}({\texttt P}_\nu)$ for its completion by the null sets of ${\texttt P}_\nu$. Moreover, write $\overline{\mathcal G} =\bigcap_{\nu} \mathcal{G}({\texttt P}_\nu)$, where the intersection is taken over all probability measures on the state space of $Y$, excluding the cemetery state.
		A finite  random time $\texttt{k}$ is called an $L$-time (generalized last exit time) if
\begin{itemize}
	\item[(i)] $\texttt{k}$ is measurable in $\overline{\mathcal G}$, and  $\texttt{k}\leq \zeta$  almost surely with respect to ${\texttt P}_\nu$, for all $\nu$,
	\item[(ii)] $\{s<\texttt{k}(\omega)-t\}=\{s<\texttt{k}(\omega_t)\}$ for all $t,s\geq 0$,
\end{itemize}	
where $\omega_t$ is the Markov shift of $\omega$ to time $t$.
The most important examples of $L$-times are killing times and last  exit times.

\begin{theorem}\label{Naga} In what follows, we work with the probability distribution
\begin{equation}
\nu(\dd a) : = \frac{\sigma_1(\dd a)|_{\mathsf{S}}}{\sigma_1(\mathsf{S})}, \qquad a \in \mathbb{R}^d,
\label{Borel}
\end{equation}
if $\mathsf{S}$ is {closed} and $\sigma_1(\mathsf{S})>0$ and, otherwise, if $\mathsf{S} = \{\vartheta\}$, $\vartheta\in \mathbb{S}^{d-1}$, we understand 
\begin{equation}
\nu(\dd a) = \delta_{\{\vartheta\} } (\dd a), \qquad a\in \mathbb{R}^d.
\label{singleton}
\end{equation}
\begin{enumerate} 
\item[(i)] For every $L$-time  $\emph{\texttt{k}}$ of $(X, \mathbb{P}^\ominus)$, the {time reversed} process
$(
X_{(\emph{\texttt{k}} - t)-}, t\ < \emph{\texttt{k}}
)$
under $\mathbb{P}_\nu^{\ominus}$ {is a time-homogeneous Markov process whose transition probabilities }
agree with those of $(X,\mathbb{P}^\vee)$.

\item[(ii)] Similarly, for every $L$-time  $\emph{\texttt{k}}$ of $(X, \mathbb{P}^\oplus)$, the {time reversed} process
$(
X_{(\emph{\texttt{k}} - t)-}, t < \emph{\texttt{k}}
)$
under $\mathbb{P}_\nu^{\oplus}$ 
is a time-homogeneous Markov process whose transition probabilities
agree with those of $(X,\mathbb{P}^\wedge)$.
\end{enumerate}
\label{main7}
\end{theorem}

{Nagasawa's result, \cite[Theorem 3.5]{Naga}, allows the definition of the time reversed process only for $t>0$, however we can extend it for $t=0$. Indeed, in (i), $\texttt{k}<\zeta=\infty$  with probability $\mathbb{P}^\ominus$ one, and the time-reversal can include $t=0$; in (ii), we may have $\texttt{k}=\zeta<\infty$  with positive probability $\mathbb{P}^\oplus$, but in this case $X_{\zeta-}=0$ with probability $\mathbb{P}^\oplus$ one, and therefore again $t=0$ can be included in the time reversed process.
	That means, if the duality is true for $t > 0$, it must be true for all $t\geq 0$.
 }

\section{Proof of Theorem \ref{main1}}

We start by recalling two useful identities. 
In Theorem 1.1 in \cite{Deep3}, the law of $X_{\underline{G}(\infty)}$ is given by
\begin{eqnarray}
\mathbb{P}_x(X_{\underline{G}(\infty)} \in \dd z) &=& c_{\alpha,d}  \frac{(|x|^2-|z|^2)^{\alpha/2}}{|z|^{\alpha}} |x-z|^{-d} \dd z, 
\qquad |x| > |z|>0,\label{closest_reach}
\end{eqnarray} 
where 
\[
c_{\alpha,d}=\pi^{-d/2}\dfrac{\Gamma\left({d}/{2}\right)^2}{\Gamma\left(({d-\alpha})/{2}\right)\Gamma\left({\alpha}/{2}\right)}.
\]
Similarly, from Corollary 1.1 of   \cite{Deep3}, it was also shown that
\begin{equation}
\mathbb{P}_x(X_{\overline{G}(\tau_1^\ominus)} \in \dd z, X_{\tau_1^\ominus} \in \dd v) = C_{\alpha,d} \frac{(|z|^2-|x|^2)^{\alpha/2}}{(|v|^2-|z|^2)^{\alpha/2} |z-v|^d |z-x|^d} \dd z \dd v,
\label{corr1.3}
\end{equation}
for $|x|<|z|<1$ and $|v|>1$, where 
\[
C_{\alpha,d}= \frac{\Gamma(d/2)^2}{\pi^{d}|\Gamma(-\alpha/2)|\Gamma(\alpha/2)}.
\]

\smallskip

First take $x \in  \bar{\mathbb{B}}^c_d$. Let  $\tau^\oplus_{\beta} := \inf \{ t>0 \colon |X_t| < \beta \}$ for any $\beta>1.$  For any $A \in \mathcal{F}_t$, define
\begin{equation}
\mathbb{P}_x^\vee (A, t < \tau^\oplus_{\beta})=\lim_{\varepsilon \rightarrow 0} \mathbb{P}_x (A, t < \tau^\oplus_{\beta} | C_\varepsilon^{\vee}).
\label{betalim}
\end{equation}
The  Markov property gives us
\begin{eqnarray}
\mathbb{P}_x (A, t < \tau^\oplus_{\beta} | C_\varepsilon^{\vee}) 
&=
\mathbb{E}_x \left[ \mathbf{1}_{\{A, t< \tau^\oplus_{\beta}\}} \dfrac{\mathbb{P}_{X_t}(C_\varepsilon^{\vee})}{\mathbb{P}_{x}(C_\varepsilon^{\vee})}\right].
\label{cond}
\end{eqnarray}
In order to prove the Theorem \ref{main1}, it is enough to prove that, for all $\beta>1,$ \eqref{out} is true for sets of the form $A \cap \{t<\tau^\oplus_{\beta}\} \in \mathcal{F}_t,$ in which case the full statement \eqref{out} follows by the Monotone Convergence Theorem as we take $\beta\downarrow1$.
Next note from \eqref{closest_reach} that
\begin{eqnarray*}
\mathbb{P}_x(X_{\underline{G}(\infty)} \in A_\varepsilon) &=& c_{\alpha,d} \int_{z \in A_\varepsilon} \frac{(|x|^2-|z|^2)^{\alpha/2}}{|z|^{\alpha}} |x-z|^{-d} \dd z \\
&=& c'_{\alpha,d} \int_1^{1+\varepsilon} \int_{\mathsf{S}} \frac{(|x|^2-r^2)^{\alpha/2}}{r^{\alpha}} |x-r\theta|^{-d} r^{d-1} \dd r\sigma_1( \dd\theta),
\end{eqnarray*}
where $c'_{\alpha,d}$ is an unimportant constant.
\smallskip

Since $(|x|^2-r^2)^{\alpha/2} |x-r\theta|^{-d}$ is continuous at $r=1$ with fixed $|x|>1$, for any $\delta>0,$ there exists $\varepsilon >0$ such that for all $1<r<1+\varepsilon$, $$(1-\delta) (|x|^2-1)^{\alpha/2} |x-\theta|^{-d}<(|x|^2-r^2)^{\alpha/2} |x-r\theta|^{-d}<(1+\delta)(|x|^2-1)^{\alpha/2} |x-\theta|^{-d}$$
and
$$\int_1^{1+\varepsilon} r^{d-\alpha+1} \dd r = c \varepsilon^{d-\alpha} +o(\varepsilon^{d-\alpha}),$$
{where $c$ is an unimportant constant.} Hence, we have 
$$\lim_{\varepsilon \rightarrow 0} \varepsilon^{\alpha-d} \mathbb{P}_x(X_{\underline{G}(\infty)} \in A_\varepsilon) = c'_{\alpha,d} \int_{\mathsf{S}} (|x|^2-1)^{\alpha/2} |x-\theta|^{-d} \sigma_1(\dd\theta),$$
where $c'_{\alpha,d}$ does not depend on $x$ {and may change from the previous one}. Note, moreover, that {for all fixed $\beta>1$}
\begin{equation}
\sup_{|x|>{\beta}}\frac{\int_{\mathsf{S}} (|x|^2-1)^{\alpha/2} |x-\theta|^{-d} \sigma_1(\dd\theta)}{|x|^{\alpha - d}}<\infty.
\label{DCT}
\end{equation}

We can both make use of the limit
\begin{equation}
\lim_{\varepsilon \rightarrow 0} \frac{\mathbb{P}_{X_t}(X_{\underline{G}(\infty)} \in A_\varepsilon)}{\mathbb{P}_{x}(X_{\underline{G}(\infty)} \in A_\varepsilon)}
= \frac{\int_{\mathsf{S}} |\theta-X_t|^{-d} (|X_t|^2-1)^{\alpha/2} \sigma_1(\dd\theta)}{\int_{\mathsf{S}} |\theta-x|^{-d} (|x|^2-1)^{\alpha/2} \sigma_1(\dd\theta)}, \qquad {t<\tau^\oplus_\beta}.
\end{equation}
as well as \eqref{DCT} {and the Dominated Convergence Theorem} to ensure the limit may be passed through the expectation in \eqref{cond} to give \eqref{out} on $\{t<\tau^\oplus_\beta\}$, thus giving the desired result.
\smallskip

Next we look at the proof of \eqref{in}. {In a similar way, it is enough to work with sets of the form $A \cap \{t<\tau^\ominus_{\beta}\} \in \mathcal{F}_t,$ with $\beta<1$.} 
 From \eqref{corr1.3}, recalling $C_\varepsilon^{\wedge}:=\{X_{\overline{G}(\tau^\ominus_1-)} \in B_{\varepsilon}\}$, we have
\begin{align}
\mathbb{P}_x(C_\varepsilon^\wedge)
&=\mathbb{P}_x (X_{\overline{G}(\tau^\ominus_1-)} \in B_{\varepsilon})\notag\\
& = C_{\alpha,d} \int_{z \in B_{\varepsilon}} \int_{v \in \mathbb{B}_d^c} \frac{(|z|^2-|x|^2)^{\alpha/2}}{(|v|^2-|z|^2)^{\alpha/2} |z-v|^d |z-x|^d} \dd z \dd v \notag\\
&=
C'_{\alpha,d} \int_{z \in B_{\varepsilon}} \frac{(|z|^2-|x|^2)^{\alpha/2}}{|z-x|^d} \dd z \int_1^{\infty} \frac{r^{d-1}\dd r}{(r^2-|z|^2)^{\alpha/2}} 
\int_{\mathbb{S}^{d-1}(0,r)}\frac{1}{|z- \theta|^{d}}\sigma_r(\dd \theta),
\label{putpoissonin}
\end{align}
where $\sigma_r(\dd \theta)$ is the surface measure on $\mathbb{S}^{d-1}(0,r)$, the sphere centred at $0$ of radius $r$, normalised to have unit mass and $C'_{\alpha,d}$ is henceforth a constant whose value may change from line to line, which depends only on $\alpha$ and $d$. The Poisson formula (giving the probability that a $d$-dimensional Brownian motion issued from $z$ (with $|z|<1$) will hit  the sphere $\mathbb{S}^{d-1}(0,r)$) tells  us that 
\begin{equation}
\int_{\mathbb{S}^{d-1}(0,r)}
 \frac{r^{d-2}(r^2-|z|^2)}{|z-\theta|^{d}}\sigma_r(\dd \theta)=1, \qquad |z|<1< r,
\label{Poisson}
\end{equation}
see for example Remark III.2.5 in \cite{KALEA}.
Putting \eqref{Poisson} in \eqref{putpoissonin} gives us
\begin{align*}
\mathbb{P}_x(C_\varepsilon^\wedge) &=  
C'_{\alpha,d} \int_{z \in B_{\varepsilon}} \frac{(|z|^2-|x|^2)^{\alpha/2}}{|z-x|^d} \dd z \int_1^{\infty} \frac{r^{d-1}}{(r^2-|z|^2)^{\alpha/2}} \frac{1}{r^{d-2}(r^2-{|z|^2})} \dd r\\
&= 
C'_{\alpha,d} \int_{z \in B_{\varepsilon}} \frac{(|z|^2-|x|^2)^{\alpha/2}}{|z-x|^d} \frac{1}{(1-|z|^2)^{\alpha/2}} \dd z \\
&= C'_{\alpha,d}\int_{1-\varepsilon}^1 \int_{\mathsf{S}} \frac{(u^2-|x|^2)^{\alpha/2}}{(1-u^2)^{\alpha/2} |u\theta-x|^d} u^{d-1} \dd u \, \sigma_1(\dd\theta  ).
\end{align*}
Since $(u^2-|x|^2)^{\alpha/2} |x-u\theta|^{-d}$ is continuous at $u=1$ with fixed $0<|x|<1$, for any $\delta>0,$ there exists $\varepsilon >0$ such that for all $1-\varepsilon<u<1$, $$(1-\delta) (1-|x|^2)^{\alpha/2} |x-\theta|^{-d}<
(u^2-|x|^2)^{\alpha/2} |x-u\theta|^{-d}<(1+\delta)(1-|x|^2)^{\alpha/2} |x-\theta|^{-d}$$
and $$\int_{1-\varepsilon}^1 \frac{u^{d-1}}{(1-u^2)^{\alpha/2}} \dd u = \int_0^\varepsilon \frac{(1-r)^{d-1}}{r^{\alpha/2}(2-r)^{\alpha/2}} \dd r = c \varepsilon^{1-\alpha/2} + o(\varepsilon^{1-\alpha/2}),$$
for an unimportant constant $c>0$.
\smallskip

It is now clear that 
$$\lim_{\varepsilon \rightarrow 0} \varepsilon^{\alpha/2-1} \mathbb{P}_x(X_{\overline{G}(\tau^\ominus_1-)} \in B_\varepsilon) = C'_{\alpha,d} \int_{\mathsf{S}} (1-|x|^2)^{\alpha/2} |x-\theta|^{-d} \sigma_1(\dd\theta).$$
Finally, we get again 
\begin{equation}
\lim_{\varepsilon \rightarrow 0} \frac{\mathbb{P}_{X_t}(X_{\overline{G}(\tau^\ominus_1-)} \in B_\varepsilon)}{\mathbb{P}_{x}(X_{\overline{G}(\tau^\ominus_1-)} \in B_\varepsilon)}
= \frac{\int_{\mathsf{S}} |\theta-X_t|^{-d} (1-|X_t|^2)^{\alpha/2} \sigma_1(\dd\theta)}{\int_{\mathsf{S}} |\theta-x|^{-d} (1-|x|^2)^{\alpha/2}\sigma_1(\dd\theta)}, \qquad {t<\tau^\ominus_\beta}.
\end{equation}
and we can proceed as in \eqref{betalim}, noting the application of dominated convergence and that for every fixed $\beta<1$,
{\begin{equation*}
	\sup_{|x|<\beta}\frac{\int_{\mathsf{S}} (|x|^2-1)^{\alpha/2} |x-\theta|^{-d} \sigma_1(\dd\theta)}{|x|^{\alpha - d}}<\infty.
	\end{equation*}}

In a similar manner, when $\mathsf{S} = \{\vartheta\},$  {we work with sets of the form $A\cap\{t<\tau^\oplus_{\beta}\}\in \mathcal{F}_t$ or $A \cap \{t<\tau^\ominus_{\beta'}\} \in \mathcal{F}_t,$ with $\beta'<1<\beta$, respectively.}  In this case, $A_\varepsilon = \{r\phi \colon r \in (1,1+\varepsilon), \phi\in\mathbb{S}^{d-1}, |\phi-\vartheta|<\varepsilon\}$ and $B_\varepsilon = \{r\phi\colon r \in (1-\varepsilon,1),  \phi\in\mathbb{S}^{d-1}, |\phi-\vartheta|<\varepsilon\}$, thus it is clear by similar analysis  that 
\begin{equation}
\lim_{\varepsilon \rightarrow 0} \frac{\mathbb{P}_{X_t}(X_{\underline{G}(\infty)} \in A_\varepsilon)}{\mathbb{P}_{x}(X_{\underline{G}(\infty)} \in A_\varepsilon)}=\lim_{\varepsilon \rightarrow 0} \frac{\mathbb{P}_{X_t}(X_{\overline{G}(\tau^\ominus_1-)} \in B_\varepsilon)}{\mathbb{P}_{x}(X_{\overline{G}(\tau^\ominus_1-)} \in B_\varepsilon)}
= \frac{|\theta-X_t|^{-d} ||X_t|^2-1|^{\alpha/2}}{ |\theta-x|^{-d} ||x|^2-1|^{\alpha/2}}.
\end{equation}
The rest of the proof is otherwise a minor adjustment of what we have seen previously, now taking account of the continuity of $(u,\theta) \mapsto |u^2-|x|^2|^{\alpha/2} |x-u\theta|^{-d}$ as well as the fact that {$\sup_{|x|>\beta} ((|x|^2-1)^{\alpha/2} |x-\theta|^{-d}  )/ |x|^{\alpha - d}<\infty$ and $\sup_{|x|<\beta'} ((1-|x|^2)^{\alpha/2} |x-\theta|^{-d}  )/ |x|^{\alpha - d}<\infty$}, in order to use dominated convergence. 
\hfill$\square$

\subsection{Proof of Proposition \ref{prop2}}

To calculate the hitting distribution, recall that  {$\mathbb{P}^{\vee}$ is the law of a stable process conditioned to attract to $\mathsf{S}$ continuously from the outside,} and $A'_\varepsilon = \{r\theta \colon r \in (1,1+\varepsilon),\theta \in \mathsf{S}'\},$ that is the restriction of $A_\varepsilon$ from the set $\mathsf{S}$ to its subset $\mathsf{S}' \subset \mathsf{S}$. Then, due to Theorem 1.3 in \cite{Deep3}, we have  {$\mathbb{P}_x^{\vee} (X_{{\texttt{k}}-} \in\mathsf{S}' )=\lim_{\varepsilon \rightarrow 0} \mathbb{P}_x (X_{\underline{G}(\infty)} \in A'_\varepsilon | C_\varepsilon^\vee)$}. Then 
\begin{eqnarray}
\lim_{\varepsilon \rightarrow 0} \mathbb{P}_x (X_{\underline{G}(\infty)} \in A'_\varepsilon | C_\varepsilon^\vee)
&=& \lim_{\varepsilon \rightarrow 0} \mathbb{P}_x (X_{\underline{G}(\infty) } \in A'_\varepsilon | X_{\underline{G}(\infty)} \in A_\varepsilon) \notag \\
&=& \lim_{\varepsilon \rightarrow 0} \frac{\mathbb{P}_x(X_{\underline{G}(\infty)} \in A'_\varepsilon)}{\mathbb{P}_x( X_{\underline{G}(\infty)} \in A_\varepsilon)}\notag \\
& =& \frac{\int_{\mathsf{S}'} |\theta-x|^{-d}\sigma_1(\dd\theta) }{\int_\mathsf{S} |\theta-x|^{-d}\sigma_1(\dd\theta)},
\end{eqnarray}
which concludes the statement in the Proposition \ref{prop2} for the case when $X$ is issued from outside. 
Similar computations give the result when 
 $X$ is issued from inside $\mathbb{B}_d$.\hfill$\square$

\section{Proof of Theorems \ref{main3} and \ref{mainx}}\label{proofmain3}

{\bf Proof of Theorem \ref{main3}:} Let us restrict our attention to the extension of $(X, \mathbb{P}^\vee)$ to include $\mathbb{S}^{d-1}\setminus\mathsf{S}$.
We need to prove that the proposed definition of $\mathbb{P}_\theta^\vee$, for any $\theta \in \mathbb{S}^{d-1} \setminus \mathsf{S}$, is is well defined as a finite entity, conforms to the correct normalisation to represent a probability measure and is consistent with the definition of $(X, \mathbb{P}^\vee)$ given in Theorem \ref{main1} on $\bar{\mathbb{B}}_d^c$, as well as offering continuous entry from the boundary $\mathbb{S}^{d-1}\setminus \mathsf{S}$.
\smallskip

We start with finiteness. To this end, we must show that, for $\theta\in \mathbb{S}^{d-1}\backslash\mathsf{S}$
\begin{equation}
\underline{\mathbb{N}}_\theta(H_\mathsf{S}(X^\epsilon(t)); t<\varsigma) <\infty, \qquad t> 0.
\label{isitfinite}
\end{equation}
Noting from \eqref{justH} that $H_{\mathsf S}(x)\leq H_{\mathbb{S}^{d-1}}(x) = |x|^{\alpha - d}(1-|x|^{-2})^{\frac{\alpha}{2}-1}$, which tends to 0 as $|x|\to \infty$, it suffices to prove that, for any $R>1$, 
\begin{equation}
\underline{\mathbb{N}}_\theta((|X^\epsilon(t)|\geq R, t<\varsigma) +\underline{\mathbb{N}}_\theta(H_\mathsf{S}(X^\epsilon(t));|X^\epsilon(t)|<R,  t<\varsigma) <\infty, \qquad t> 0.
\label{isitfinite2}
\end{equation}
Abusing notation and using $\tau^\ominus_R = \inf\{t>0: |X^\epsilon_t|>R\}$ in the canonical sense,
\begin{align}
\underline{\mathbb{N}}_\theta((|X^\epsilon(t)|\geq R, t<\varsigma) &\leq \underline{\mathbb{N}}_\theta(t<\varsigma\wedge \tau^\ominus_R) \leq \underline{n} (t< \kappa_{\log R}\wedge\zeta)<\infty
\label{bigexcursion}
\end{align}
where $\kappa_{\log R}= \inf\{t>0: \epsilon(t)>\log R\}$ and $\underline{n}$ is the excursion measure of $\xi-\underline\xi$. (The fact that the second expression in \eqref{bigexcursion} is finite is a well known fact from the theory of L\'evy processes;  otherwise there would be an infinite rate of having arbitrary large excursions, which occurs with probability zero.)
\smallskip

Our objective now will be to show that in fact the resolvent 
\begin{equation}
I_\textsf{S}(\theta): = \int_0^\infty \underline{\mathbb{N}}_\theta(H_\mathsf{S}(X^\epsilon(t));|X^\epsilon(t)|<R,  t<\varsigma) \dd t<\infty,
\label{specialresolvent}
\end{equation}
which ensures that \eqref{isitfinite} is finite for Lebesgue almost all $t>0$ and hence, thanks to stochastic continuity of the excursion measure, for all $t>0$.

To prove \eqref{specialresolvent} consider $|y|\geq 1$ and $\theta\in \mathbb{S}^{d-1} \setminus \mathsf{S}$, we can appeal to Proposition 5.2 of \cite{Deep3}, 
which identifies, 
for $x\in\mathbb{R}^d\backslash\{0\}$, and continuous  $g:\mathbb{R}^d\mapsto\mathbb{R}$ whose support is compactly embedded in the exterior of the ball of radius $|x|$,
\begin{align}
\underline{\mathbb{N}}_{\arg(x)}\left(\int_0^{\zeta}g(|x|{\rm e}^{\epsilon(u)}\Theta^\epsilon(u)){\rm d} u\right)
&= 
\frac{\Gamma((d - \alpha)/2)}{2^{\alpha} \pi^{d/2}\Gamma(\alpha/2)}
\int_{|x|<|z|} g(z) 
\frac{(|z|^2-|x|^2)^{\alpha/2}}
{|z|^\alpha|x-z|^{d} }\dd z.
\label{upNint}
\end{align}
This gives us
\begin{align*}
\rho^\vee(\theta, \dd y)& : = \int_0^\infty \mathbb{P}_\theta^\vee(X_t\in \dd y, t<\mathtt{k})\dd t\\
& =\frac{\Gamma(d/2)}{\Gamma((d-\alpha)/2)\Gamma(\alpha/2+1)}\frac{H_\mathsf{S}( y)}{h(\theta)} \int_0^\infty  \underline{\mathbb{N}}_\theta \left(X^\epsilon(t)\in {\rm d}y, \, t<\varsigma \right)\dd t\\
&
=\frac{\Gamma(d/2)}{\Gamma((d-\alpha)/2)\Gamma(\alpha/2+1)}\frac{H_\mathsf{S}( y)}{h(\theta)} |y|^{\alpha}\int_0^\infty  \underline{\mathbb{N}}_\theta \left({\rm e}^{\epsilon(u)}\Theta^\epsilon(u)\in {\rm d}y, \, u<\zeta \right)\dd u
\\
&= 
\frac{\Gamma(d/2)}{2^{\alpha} \pi^{d/2}\Gamma(\alpha/2)\Gamma(\alpha/2+1)}
\frac{||y|^2-1|^{\alpha/2}}{|\theta-y|^{d} }\frac{H_\mathsf{S}( y)}{h(\theta)}\dd y\\
&=\frac{\Gamma(d/2)}{2^{\alpha} \pi^{d/2}\Gamma(\alpha/2)\Gamma(\alpha/2+1)}
\frac{H_\mathsf{\{\theta\}}( y)H_\mathsf{S}( y)}{h(\theta)}\dd y
\end{align*}
where  we recall $h(\theta) = \int_{\mathsf{S}}|\theta -\vartheta|^{-d}\sigma_1(\dd \vartheta)$, the representation of $X^\epsilon$ is given in  \eqref{Xexcursions} and the fact that $ {\rm e}^{\alpha\epsilon(\varphi(t))} \dd \varphi(t) =|X^\epsilon_t|^{\alpha}  \dd \varphi(t)= \dd t$ on $t<\varsigma$.
\smallskip

It now follows that, up to a multiplicative constant $C$ (which in the following calculations will play the role of different constants that may change from line to line)
\begin{align}
I_\textsf{S}(\theta) &= \int_{1<|y|<R}\rho^\vee(\theta, \dd y)\notag\\
&=\frac{C}{h(\theta)}\int_{1<|y|<R}H_\mathsf{\{\theta\}}( y)H_\mathsf{S}( y)\dd y\notag\\
&=\frac{C}{h(\theta)}
 \int_{1<|y|<R} \int_{\phi \in \mathsf{S}} \frac{||y|^2-1|^\alpha}{|\theta-y|^d| \phi-y|^d} \dd y \sigma_1(\dd\phi).
 \label{IS}
\end{align}
Since $ \mathbb{S}^{d-1}\backslash\mathsf{S}$ is open, it is easy to see that we can choose $\varepsilon$ small enough such that, for $\theta\in\mathbb{S}^{d-1}\backslash\mathsf{S}$, 
 \begin{align}
 \{y\in\mathbb{R}^{d}: |y|>1\}& = \{y\in\mathbb{R}^{d}:|y|>1\text{ and } |y-\theta|>\varepsilon\}
 \notag\notag\\
&\hspace{2cm} \cup \{y\in\mathbb{R}^{d}: |y|>1\text{ and }|y-\phi|>\varepsilon, \text{ for all }\phi\in\mathsf{S}\},
\label{split}
 \end{align}
 such that 
 \begin{align*}
&  \{y\in\mathbb{R}^{d}:|y|>1\text{ and } |y-\theta|\leq \varepsilon\} \\
&\hspace{2cm}  \cap\{y\in\mathbb{R}^{d}:|y|>1\text{ and } |y-\phi|\leq \varepsilon, \text{ for some }\phi\in\mathsf{S}\} = \emptyset.
 \end{align*}
Making use of \eqref{lessthanone}, \eqref{justH} and \eqref{3.665}, and that, for $r>1$, ${_2}F_1(d/2, 1; d/2, r^{-2}) = (1-r^{-2})^{-1}$, allowing $C$ to again play the role of a strictly positive constant that may change from line to line, we have, for $\theta\not\in \mathsf{S}$,
\begin{eqnarray*}
I_{\mathsf{S}}(\theta) &=& \frac{C}{h(\theta)}\int_{1<|y|<R} \int_{\phi \in \mathsf{S}} \frac{||y|^2-1|^\alpha}{|\theta-y|^d| \phi-y|^d} \dd y \sigma_1(\dd \phi)\\
&\leq& \frac{C}{h(\theta)}\int_{1<|y|<R,|\theta-y| \geq \varepsilon} \int_{\phi \in \mathsf{S}} \frac{||y|^{2}-1|^\alpha}{|\theta-y|^d| \phi-y|^d} \dd y \sigma_1(\dd \phi)\\
&&\hspace{2cm}+\frac{C}{h(\theta)} \int_{1<|y|<R,|\phi-y| \geq \varepsilon} \int_{\phi \in \mathsf{S}} \frac{||y|^{2}-1|^\alpha}{|\theta-y|^d| \phi-y|^d} \dd y \sigma_1(\dd \phi)\\
&\le& \varepsilon^{-d} \frac{C}{h(\theta)}\Big(\int_{1<|y|<R} \int_{\phi \in \mathsf{S}} \frac{||y|^{2}-1|^\alpha}{|\phi-y|^d} \dd y \sigma_1(\dd \phi)+ \int_{1<|y|<R} \int_{\phi \in \mathsf{S}} \frac{||y|^{2}-1|^\alpha}{|\theta-y|^d} \dd y\sigma_1( \dd \phi)\Big) \\
&\le& \varepsilon^{-d} \frac{C}{h(\theta)}\Big(\int_{1<|y|<R} \int_{\phi \in \mathbb{S}^{d-1}} \frac{||y|^{2}-1|^\alpha}{|\phi-y|^d} \dd y \sigma_1(\dd \phi)+ \sigma_1 (\mathsf{S}) \int_{1<|y|<R}  \frac{||y|^{2}-1|^\alpha}{|\theta-y|^d} \dd y \Big) \\
&=& \varepsilon^{-d} \frac{C}{h(\theta)} \Big(  \int_{1<|y|<R}  ||y|^2-1|^\alpha |y|^{\alpha-d} (1-|y|^{-2})^{\frac{\alpha}{2}-1}\dd y  \notag\\
&&\hspace{2cm}+ \int_{1}^R  \int_0^\pi \frac{r^{d-1}(r^2-1)^\alpha (\sin\vartheta)^{d-2}}{(r^2-2r\cos(\vartheta)+1)^{d/2}} \dd r \dd\vartheta\Big) \\
&=& \varepsilon^{-d} \frac{C}{h(\theta)} \left(\int_{1}^R (r^2-1)^{\frac{3\alpha}{2}-1} 
r\dd r
+\int_{1}^R (r^2-1)^{\alpha-1} r\dd r \right) \\
&=&   \varepsilon^{-d} \frac{C}{h(\theta)} \left(\int_{1}^{R^2}  (u-1)^{\frac{3\alpha}{2}-1}\dd u  + \int_{1}^{R^2}  (u-1)^{\alpha-1}\dd u\right)<\infty.
\end{eqnarray*}
\smallskip

{Now let us turn to the issue of consistency. Recall that $(\Lambda, \underline{\mathbb{N}}_{\cdot})$ is an exit system for the process $(\xi,\Theta)$. In particular, under any measure $\underline{\mathbb{N}}_{\theta}$ the process $((\epsilon(s), \Theta^\epsilon(s)), s<\zeta)$ is  a strong Markov process with the same semigroup as $(\xi, \Theta)$ killed at its first hitting time of $(-\infty,0]\times\mathbb{S}^{d-1}$, see \cite[Theorem 6.3]{Maison}.}
As a consequence,  for $\theta\in\mathbb{S}^{d-1}\setminus \mathsf{S}$, 
\begin{eqnarray}
\mathbb{E}_\theta^\vee[g(X_{t+s})] &=& \frac{C}{h(\theta)}\underline{\mathbb{N}}_\theta (H_\mathsf{S}(X_{t+s}^\epsilon) g(X^\epsilon_{t+s}) \mathbf{1}_{(s+t < \varsigma)}) \notag \\
&=& \frac{C}{h(\theta)}\underline{\mathbb{N}}_\theta \left(H_\mathsf{S}(X_t^\epsilon) \mathbf{1}_{(t < \varsigma)} \mathbb{E}_{X_t^\epsilon} \left [ \frac{H_\mathsf{S}(X_{s})}{H_\mathsf{S}(X_t^\epsilon)} g(X_{s}) \mathbf{1}_{(s <\tau^\ominus_1)} \right] \right)\notag \\
&=&\frac{C}{h(\theta)} \underline{\mathbb{N}}_\theta \left(H_\mathsf{S}(X_t^\epsilon) \mathbf{1}_{(t < \varsigma)} \mathbb{E}_{X_t^\epsilon}^\vee[g(X_s)] \right)\notag\\
&=&\mathbb{E}_\theta^\vee\left[H_\mathsf{S}(X_t^\epsilon) \mathbf{1}_{(t < \varsigma)} \mathbb{E}_{X_t^\epsilon}^\vee[g(X_s)] \right],
\end{eqnarray} 
where $C =  {\Gamma(d/2) } /{\Gamma(\alpha/2+1)\Gamma((d-\alpha)/2) }$.
Hence, using the notation $\mathcal{P}_t^\vee[g](x):= \mathbb{E}_x^\vee [g(X_t)],$ we have $\mathcal{P}_{t+s}^\vee[g](x)=\mathcal{P}_{t}^\vee [\mathcal{P}_{s}^\vee [g]](x)$ for any $x \in \mathbb{R}^d \setminus (\mathbb{B}_d \cup \mathsf{S}),$ and the required consistency follows.  
\smallskip

 To demonstrate the consistent choice of   normalisation in our definition of $\mathbb{P}^\vee$, we will reconsider a different derivation of the resolvent $\rho^\vee$. To this end, suppose that $x\in\bar{\mathbb{B}}^c$ and we can similarly consider the resolvent of $(X,\mathbb{P}_x^\vee)$.  This calculation can be developed using  the nature of the Doob $h$-transform \eqref{out} and Theorem III.3.4 in \cite{KALEA} and takes the form
\begin{align}
\label{resolventawayfromsphere}
\rho^\vee(x,\dd y) = \frac{H_{\mathsf{S}}(y)}{H_\mathsf{S}(x)}\rho^\oplus(x,\dd y),\qquad |x|,|y|>1,
\end{align}
where 
\begin{equation}
\label{rhoplus}
\rho^\oplus(x,\dd y)=\frac{\Gamma(d/2)}{2^{\alpha}\pi^{d/2}\Gamma(\alpha/2)^2}|x-y|^{\alpha - d} \int_0^{\zeta^\oplus(x,y)}  (u+1)^{-d/2}u^{\alpha/2-1}\dd u \, \dd y
\end{equation}
and $\zeta^\oplus(x,y) = (|x|^2-1)(|y|^2-1)/|x-y|^{2}$. To show continuity as $x\to\theta\in \mathbb{S}^{d-1}\backslash\mathsf{S}$, and hence that the choice of normalisation in \eqref{veeS} is correct,  we note that, as $r\to1$, 
\begin{align}
\frac{H_\mathsf{S}(y)}{H_{\mathsf{S}}(x)}\rho^\oplus(r\theta, \dd y)&\sim\frac{\Gamma(d/2)|r\theta - y|^{\alpha -d}H_\mathsf{S}(y)}{2^{\alpha}\pi^{d/2}\Gamma(\alpha/2)^2
h(\theta)}\notag\\
&\hspace{1cm}\times\frac{ 2r (|y|^2-1)  |r\theta-y|^{-2}\zeta(r\theta,y)^{\alpha/2-1} (1+\zeta(r\theta,y))^{-d/2}}{2r(\alpha/2)(r^2-1)^{\alpha/2-1}}
\, \dd y\notag\\
&\sim \frac{\Gamma(d/2)
(|y|^2-1)^{\alpha/2}|\theta - y|^{-d}
H_\mathsf{S}(y)}{2^{\alpha}\pi^{d/2}\Gamma(\alpha/2)\Gamma(\alpha/2 + 1)
h(\theta)}\dd y\notag\\
&=\frac{\Gamma(d/2)}{2^{\alpha} \pi^{d/2}\Gamma(\alpha/2)\Gamma(\alpha/2+1)}
\frac{H_\mathsf{\{\theta\}}( y)H_\mathsf{S}( y)}{h(\theta)}\dd y
\end{align}

Now, we need to show that $\mathbb{P}_\theta^\vee (X_{0+} = \theta) =1$ for any $\theta \in \mathbb{S}^{d-1} \setminus \mathsf{S}.$ Since $\lim_{t \downarrow 0} \varphi(t) = 0,$ it suffices to show that
\begin{equation}
    \mathbb{P}_\theta^\vee (X_0 \ne \theta) = \underline{\mathbb{N}}_\theta \left(\textstyle{\left\{\lim_{t \downarrow 0} \epsilon(t) = 0, \,\lim_{t \downarrow 0} \Theta^\epsilon (t) = \theta \right\}^c}\right) = 0.
    \label{zero}
\end{equation}
Let us first observe $\epsilon$ is an excursion of $\xi$ from its running minimum and $\xi$ is a hypergeometric L\'evy process with unbounded variation, hence $0$ is regular for $(0,\infty)$, that is  
\begin{equation*}
    \mathbf{P}_{0,\theta}(\tau_0^+ = 0) =1, \qquad \theta\in \mathbb{S}^{d-1},
\end{equation*}
where $\tau_0^+ = \inf\{t>0 : \xi_t > 0\}.$ Classical  excursion theory for L\'evy processes implies that 
the excursions of $\xi$ from its infimum begin continuously. Thanks to isotropy, this is equivalent to saying
\begin{equation}
\underline{\mathbb{N}}_\theta \left(\textstyle{\{\lim_{t \downarrow 0} \epsilon(t) = 0\}^c}\right)=0.
\label{eps}
\end{equation}
Since the jump measure of $X$ in radial form is
\begin{equation*}
    \Pi (\dd r,\dd\theta) = \frac{1}{r^{1+\alpha}} \sigma_1 (\dd\theta){\rm d}r ,\qquad r>0, \theta\in\mathbb{S}^{d-1},
\end{equation*} 
as a consequence, the process $(\xi, \Theta)$ has the property that both the modulator and the ordinate must jump simultaneously (the precise jump rate was explored in \cite{KALEA}). If it were the case that $\underline{\mathbb{N}}_\theta \left(\textstyle{\left\{\lim_{t \downarrow 0} \Theta^\epsilon (t) = \theta \right\}^c}\right)>0$ (and hence for all $\theta\in\mathbb{S}^{d-1}$ by rotational symmetry), this would be tantamount to a discontinuity in $\Theta$ but not in $\xi$, which is a contradiction {since  $((\epsilon(s), \Theta^\epsilon(s)), s<\zeta)$  under $\underline{\mathbb{N}}_{\theta}$ has the same semigroup as the isotropic process $(\xi, \Theta)$ killed at its first hitting time of $(-\infty,0]\times\mathbb{S}^{d-1}.$} The requirement \eqref{zero} now follows. This completes the proof of Theorem \ref{main3} as far as $\mathbb{P}^\vee$ is concerned. 

\smallskip

The proof of Theorem \ref{main3} for $(X, \mathbb{P}^\wedge)$  is essentially the same as soon as we have an analogous identity for \eqref{upNint},  but for $\overline{\mathbb N}_\theta$. Unfortunately this does not seem to be available in the literature, and so we spend a little time developing it here. However the remaining details of the proof of Theorem \ref{main3} we leave to the reader. 

The main idea behind the derivation of an analogue to \eqref{upNint} for $\overline{\mathbb N}_\theta$ lies with the use of the Riesz--Bogdan-\.Zak transform in Theorem \ref{RBSthrm}. Let us consider a variant of the radial excursion process which is based on the MAP $(\overline\xi - \xi, \Theta)$, that is associated to $X$ but now under the change of measure \eqref{COM}. The reader will recall that the probabilities $\mathbb{P}^\circ =(\mathbb{P}^\circ_x, x\neq 0)$ correspond to conditioning the process $X$ to be continuously absorbed at the origin. It turns out that  if $(\xi, \Theta)$, with probabilities $\mathbf{P}^\circ =(\mathbf{P}^\circ_{x,\theta}, x\in\mathbb{R}, \theta\in\mathbb{S}^{d-1})$, is the MAP whose Lamperti transform gives  $(X, \mathbb{P}^\circ)$, then  $(-\xi,\Theta)$, is the MAP whose Lamperti transform gives $(X, \mathbb{P})$;
see Theorem 1.3.6 and Corollary 1.3.17 of \cite{KALEA}.

\smallskip

In the spirit of \eqref{exitsystem1} we can write down the exit system for the radial excursion process of $(\overline{\xi}- \xi, \Theta)$ from $\{0\}\times\mathbb{S}^{d-1}$ under $\mathbb{P}^\circ$. Suppose that $\ell^\circ$, and  $(\overline{\mathbb N}^\circ_\theta, \theta\in\mathbb{S}^{d-1}))$ denote the associated local time and system of excursion measures. As with  excursion theory from the radial minimum of $X$, isotropy  allows us to conclude that we may  choose $\ell^\circ$ to be the local time at 0 of $\overline\xi - \xi$, and that $\xi$ (without its modulator $\Theta$) is necessarily a L\'evy process under $\mathbf{P}$. Since  $\lim_{t\to\infty}\xi_t = -\infty$ under $\mathbf{P}^\circ$, we can also appeal to isotropy again to normalise $\ell^\circ$ in such a way that $\overline{\mathbb N}^\circ(\zeta = \infty) = 1$.

\smallskip

With this set up we can follow the reasoning in \cite{Deep3} and deduce that, for positive, bounded and measurable $g$ on $\mathbb{R}^d$,
\begin{equation}
\label{Ncirc}
\overline{\mathbb{N}}^\circ_{\arg(x)}\left(\int_0^\zeta g({\rm e}^{\epsilon(s)} \Theta^\epsilon(s)) {\rm e}^{\alpha\epsilon(s)} \dd s\right)
= \lim_{|x|\uparrow1} \frac{
\mathbb{E}^\circ_x\left(\int_0^{\tau_1^\ominus} g(X_s)\dd s\right)
}{
\mathbb{P}^\circ_x(\tau^\ominus_1 = \infty)
},
\end{equation}
where we recall that $\tau^\ominus_1 = \inf\{t>0: |X_t|>1\}$. Note that the choice of normalisation of $\ell^\circ$ is implicit in the aforementioned limiting equality. Appealing to numerous calculations involving the Riesz--Bogdan--\.Zak transformation e.g. in \cite{KALEA} or indeed \cite{AEKJC}, we can rewrite the limit 
\[
\lim_{|x|\uparrow1} \frac{
\mathbb{E}^\circ_x\left(\int_0^{\tau_1^\ominus} g(X_s)\dd s\right)
}{
\mathbb{P}^\circ_x(\tau^\ominus_1 = \infty)
}
= \lim_{|x|\uparrow1} \frac{
\mathbb{E}_{Kx}\left(\int_0^{\tau_1^\oplus} g(KX_s)|X_s|^{-2\alpha}\dd s\right)
}{
\mathbb{P}_{Kx}(\tau^\oplus_1 = \infty),
}
\]
where $Kx = x/|x|^2$.
Appealing to the identities provided in \eqref{z=oo} and \eqref{rhoplus}, the limiting ratio is computable directly giving us in \eqref{Ncirc}
\[
\overline{\mathbb{N}}^\circ_{\arg(x)}\left(\int_0^\zeta g({\rm e}^{\epsilon(s)} \Theta^\epsilon(s)) {\rm e}^{\alpha\epsilon(s)} \dd s\right)
= \int_{|z|>1} g(Kz) |z|^{-2\alpha} \frac{(|z|^2-1)^{\alpha/2} }{|\arg(x) - z|^d}\dd z.
\]
An easy change of variables $y = Kz$, noting the classical analytical facts that $\dd z = |y|^{-2d}\dd y$ and $|\theta - Ky| = |\theta - y|/|y|$, for $\theta\in\mathbb{S}^{d-1}$, 
\begin{equation}
\overline{\mathbb{N}}^\circ_{\theta}\left(\int_0^\zeta g({\rm e}^{\epsilon(s)} \Theta^\epsilon(s)) {\rm e}^{\alpha\epsilon(s)} \dd s\right)
= \int_{|y|<1} g(y) |y|^{\alpha-d} \frac{(1-|y|^2)^{\alpha/2} }{|\theta - y|^d}\dd z.
\label{forcomparison}
\end{equation}

As noted in \cite{KALEA}, the change of measure \eqref{COM} 
when understood as a change of measure affecting $(\xi, \Theta)$, is equivalent to the martingale change of measure,
\begin{equation}
\label{boldCOM}
\left.\frac{\dd \mathbf{P}_{x,\theta}^\circ}{\dd \mathbf{P}_{x,\theta}}\right|_{\sigma((\xi_s, \Theta_s), s\leq t)} = {\rm e}^{(\alpha - d)(\xi_t-x)}.
\end{equation}
We can use this to compare the left-hand side of \eqref{forcomparison} with an analogous object albeit for $\overline{\mathbb N}_\theta$, the excursion measure of $(\overline{\xi}-\xi, \Theta)$ from $\{0\}\times \mathbb{S}^{d-1}$ under $\mathbf{P}$, by studying the effect of \eqref{boldCOM} on the exit formula \eqref{exitsystem1}. It is straightforward to show that, for $\theta\in\mathbb{S}^{d-1}$ and positive, bounded and measurable $g$,
\[
\overline{\mathbb{N}}^\circ_\theta\left(\int_0^\zeta g({\rm e}^{\epsilon(s)} \Theta^\epsilon(s)) \dd s\right) = 
\overline{\mathbb{N}}_\theta\left(\int_0^\zeta g({\rm e}^{\epsilon(s)} \Theta^\epsilon(s)) {\rm e}^{(\alpha -d)\epsilon(s)}\dd s\right) .
\]
Note that the normalisation of local time for $(\xi,\Theta)$ under $\mathbf{P}$ is, in effect, chosen by the above equality. It follows that 
\begin{equation}
\label{barN}
\overline{\mathbb{N}}_\theta\left(\int_0^\zeta g({\rm e}^{\epsilon(s)} \Theta^\epsilon(s)) \dd s\right)= \int_{|y|<1} g(y) \frac{(1-|y|^2)^{\alpha/2}}{|y|^\alpha |\theta - y|^d}\dd y.
\end{equation}
The reader will note that, aside from the domain of integration on the right-hand side, this agrees with \eqref{upNint}. 

\smallskip

With \eqref{barN} in hand, as alluded to above, we can now leave the reader to verify that the proof of  Theorem \ref{main3} for $(X, \mathbb{P}^\wedge)$ is essentially verbatim the same as for $(X, \mathbb{P}^\vee)$.
\hfill $\square$

\bigskip
\noindent{\bf Proof of Theorem \ref{mainx}:}
Given the proof of Theorem \ref{main3} above, we refrain from giving the proof of Theorem \ref{mainx}, noting only that it is a variant of the arguments given there. The details are, once again,  left to the reader.  We additionally note that e.g.  in this case of $\mathbb{P}^\ominus$, the excursion may begin anywhere on $\mathbb{S}^{d-1}$ and, when proving that e.g. $\underline{\mathbb{N}}_\theta(H^\ominus(X^\epsilon(t)); t<\varsigma) <\infty$,  it is much easier to show that the analogue of \eqref{specialresolvent} is finite without needing to split  space up as in \eqref{split}.\hfill$\square$

\section{Proof of Theorem \ref{main7}}
Recall the notation for a general Markov process $(Y,\texttt{P})$ on $E$ preceding the statement of Theorem \ref{Naga}. 
We will additionally write $\mathcal{P}: = (\mathcal{P}_t, t\geq 0)$ for the   semigroup associated to $(Y,\texttt{P})$.

\smallskip

 Theorem 3.5 of Nagasawa \cite{Naga}, shows that, under suitable assumptions on the Markov process, $L$-times form a natural family of random times  at which the pathwise time-reversal 
\[
\stackrel{_\leftarrow}{Y}_{\!t}:=Y_{(\texttt{k}-t)-},\qquad  t\in {(0,\texttt{k})},
\]
 is again a Markov process.  Let us state Nagasawa's principle assumptions.

\smallskip

\textbf{(A)} The potential measure $U_Y(a, \cdot)$ associated to $\mathcal{P}$, defined by the relation 
\begin{equation}
\int_Ef(x)U_Y(a,\dd x) = \int_0^\infty \mathcal{P}_t[f](a)\dd t={\texttt E}_a\left[\int_0^\infty f(X_t)\,\dd t\right],\qquad a\in E,
\label{GY}
\end{equation}
 for bounded and measurable $f$ on $E$, is $\sigma$-finite. Assume  that there exists  a probability measure, $\nu$, such that, if we put
\begin{align}\label{a1}
	\mu(A)=\int U_Y(a,A)\, \nu(\dd a)\quad \text{ for }A\in \mathcal B(\R),
\end{align}
	then there exists a Markov transition semigroup, say $\hat{\mathcal{P}}: = (\hat{\mathcal P}_t, t\geq 0)$ such that 	\begin{align}
		\int_E \mathcal{P}_t[f](x) g(x)\, \mu(\dd x)=\int_E f(x) \hat{\mathcal P}_t [g](x)\,\mu(\dd x),\quad t\geq 0,
		\label{weakdualtity}
	\end{align}
	for bounded, measurable and compactly supported test-functions $f, g$ on $E$.\smallskip
	


\smallskip

\textbf{(B)} For any continuous test-function $f\in C_0(E)$, the space of continuous and compactly supported functions,  and $a\in E$, assume that $\mathcal{P}_t[f](a)$ is right-continuous in $t$ for all $a\in E$ and, for $q> 0$, ${U}_{\hat Y}^{(q)}[f](\stackrel{_\leftarrow}{Y}_{\!t})$ is right-continuous in $t$, where, for bounded and measurable $f$ on $E$,
 \[
{U}_{\hat Y}^{(q)}[f](a) =\int_0^\infty  {\rm e}^{-qt}\hat{\mathcal{P}}_t[f](a)d t,\qquad  a\in E,\]
is the $q$-potential associated to $\hat{\mathcal P}$.
\smallskip

 Nagasawa's  duality theorem, Theorem 3.5. of \cite{Naga}, now reads as follows.

 \begin{theorem}[Nagasawa's duality theorem]\label{Ndual} Suppose that assumptions {\rm{\bf (A)} } and {\rm{\bf (B)}} hold. For the given starting probability distribution $\nu$ in {\rm{\bf (A)} } and any $L$-time $\emph{\texttt{k}}$, the time-reversed process $\stackrel{_\leftarrow}{Y}$ under $\emph{\texttt P}_\nu$ is a time-homogeneous Markov process with transition probabilities
\begin{align}
	\emph{\texttt{P}}_\nu(\stackrel{_\leftarrow}{Y}_t \in A\,|\stackrel{_\leftarrow}{Y}_r, 0<r< s)=\emph{\texttt{P}}_\nu(\stackrel{_\leftarrow}{Y}_t \in A\,|\stackrel{_\leftarrow}{Y}_s)={p}_{\hat{Y}}(t-s,\stackrel{_\leftarrow}{Y}_s,A),\quad \emph{\texttt{P}}_\nu\text{-almost surely},
\end{align}
for all $0<s<t$ and Borel $A$ in $\mathbb{R}$, where ${p}_{\hat{Y}}(u, x, A)$, $u\geq 0$, $x\in\mathbb{R}$, is the transition measure associated to the semigroup $\hat{\mathcal P}$.
\end{theorem}

\noindent {\bf Proof of Theorem \ref{main7}}. 
We give the proof of (i), the proof of (ii) is almost identical albeit requiring some straightforward adjustments. Once again, we leave the details to the reader. {When $t > 0$, we use Nagasawa's duality theorem.    However, since the process is conditioned to hit 
	continuously, its dual processes from the hitting time must leave the sphere continuously. That
means, if the duality is true for $t > 0$, it must be true for all $t\geq 0$.}

\smallskip

We will make a direct application of Theorem \ref{Ndual}, with $Y$ taken to be the process $(X,\mathbb{P}^\ominus_\nu)$ where $\nu$ satisfies \eqref{Borel} or \eqref{singleton} according to the nature of $\mathsf{S}$. Accordingly, we will write $U^\ominus$ in place of $U_Y$, $\mathcal{P}^\ominus$ in place of $\mathcal{P}$ etc.  Moreover, the dual process, formerly $\hat Y$,  is taken to be $(X,\mathbb{P}^\vee)$ and we will, in the obvious way, work with the notation $U^\vee$ in place of $U_{\hat{Y}}$, $\mathcal{P}^\vee$ in place of $\hat{\mathcal{P}}$ and so on. In essence we need only to verify the two assumptions {\bf (A)} and {\bf (B)}.
Let us momentarily take the former of these two cases. 
\smallskip

In order to verify {\bf (A)} we will make use of \eqref{upNint}. 
Noting that  ${\rm e}^{\alpha \epsilon_{\varphi(t)}} \dd\varphi(t)=\dd t$, we have  for $a\in \mathbb{S}^{d-1}\setminus \mathsf{S}$ and bounded measurable $f:\mathbb{R}^d\setminus (\mathbb{B}_d\cup \mathsf{S})\to [0,\infty)$,
\begin{align}
 U^\ominus[f](a) &= \mathbb{E}^\ominus_a\left[\int_0^\infty f(X_t ) \dd t \right] \notag\\
    &= \underline{\mathbb{N}}_a \left(\int_0^\varsigma  H^\ominus(X^\epsilon_t) f(X^\epsilon_t)\dd t \right) \notag\\
    &= \underline{ \mathbb{N}}_a \left(\int_0^\varsigma 
    H^\ominus({\rm e}^{\epsilon(u)} \Theta^\epsilon(u) ) f({\rm e}^{\epsilon(u)} \Theta^\epsilon(u) )  {\rm e}^{\alpha \epsilon_{u}} \dd u \right)\notag \\
    &=C\int_{\mathbb{R}^d\setminus (\mathbb{B}_d\cup \mathsf{S})}H^\ominus(y) f(y) (|y|^2-1)^{\alpha/2} |a-y|^{-d} \dd y,
    \label{Gominus}
\end{align}
where $U^\ominus[f](a) = \int_{\mathbb{R}^d\setminus (\mathbb{B}_d\cup \mathsf{S})} f(y)U^\ominus(a,\dd y)$, $C>0$ is an unimportant constant and we have used  \eqref{Nominus} in the second equality.
\smallskip

Next, we need to develop an expression for the reference measure $\mu$. This only needs to  be identified  up to a multiplicative constant. As such, in the setting that  $\sigma_1(\mathsf{S})>0$, recalling \eqref{a1}, \eqref{Borel} and \eqref{H_S}, we can take (ignoring multiplicative constants in each line)
\begin{align}
   \mu(\dd y) &=\int_{\mathsf{S}}\nu(\dd a) U^\ominus(a,\dd y) \notag\\
   &=\int_{\mathsf{S}}\sigma_1(\dd a) H^\ominus(y)  (|y|^2-1)^{\alpha/2} |a-y|^{-d} \dd y\notag\\
&   =H_{\mathsf{S}}(y) H^\ominus(y) \dd y, \qquad y\in \mathbb{R}^d\setminus (\mathbb{B}_d\cup \mathsf{S}).
\end{align}
When $\mathsf{S} = \{\vartheta\}$, we replace the use of \eqref{Borel} by \eqref{singleton} in the above calculation and the same answer comes out (up to a multiplicative constant).
\smallskip

Next, we need to verify that  \eqref{weakdualtity} holds. Indeed, using Hunt's switching identity (cf. Chapter II.1 of  \cite{Bertoin}) for the  process $(X_t , t<\tau^\oplus_1)$, we have  for $x,y\in\mathbb{R}^d\setminus \bar{\mathbb{B}}_d$
\begin{align*}
 \mu(\dd y) \mathcal{P}_t^{\ominus} (y,\dd x) 
      &=  \mathcal{P}_t^{\ominus} (y,\dd x) H_\mathsf{S}(y) H^\ominus(y)\dd y  \\
            &=  \frac{H^\ominus(x)}{H^\ominus(y)}\mathcal{P}_t^{\mathbb{B}_d} (y,\dd x)  H_\mathsf{S}(y) H^\ominus(y)\dd y  \\
    &=  \mathcal{P}_t^{\mathbb{B}_d} (x,\dd y) H_\mathsf{S}(y) H^\ominus(x) \dd x \\
    &=    \mathcal{P}_t^\vee (x,\dd y) \mu(\dd x),
\end{align*}
where $\mathcal{P}_t^{\mathbb{B}_d} (x,\dd y)  = \mathbb{P}_x(X_t \in \dd y, \, t<\tau^\oplus_1)$.
Note, as the measure $\mu$ is absolutely continuous with respect to Lebesgue measure, we do not need to deal with the case that $x$ or $y$ belong to $\mathbb{S}^{d-1}\setminus\mathsf{S}$.
\smallskip

Let us now turn to the verification of assumption {\bf (B)}.
This assumption is immediately satisfied on account of the fact that both $\mathcal{P}^\ominus$ and $\mathcal{P}^\vee$ are right-continuous semigroups by virtue of their definition as a Doob $h$-transform with respect to the Feller semigroup $\mathcal{P}^{\mathbb{B}_d} $ of the stable process killed on entry to $\mathbb{B}_d$. 
With both {\bf (A)} and {\bf (B)} in hand, we can invoke Theorem \ref{Ndual} and the desired result follows.\hfill$\square$
%
%
%

\section{Concluding remarks}
The results in this paper have considered the setting of conditioning a relatively special class of Markov process to  continuously hit a subset of the unit sphere with a one-sided approach. Taking a step back, one would ideally like to drop a number of the specialisms specific to our approach e.g. moving to a general Markov process and conditioning it continuously hit   a suitably general domain.  The current proofs rely on too many particular features of stable L\'evy processes for the results to directly generalise in this respect. For example, suppose that we drop the assumption that the stable process continuously approaches $\mathsf{S}$ from just one side, but instead we allow it to continuously approach without radial confinement. This is a topic that has been addressed in follow-on work \cite{Tsogii2}, for which  a mixture of features that are specific to stable L\'evy processes together with  general potential-analytic considerations are used. The classical work of Doob \cite{doob} for the setting of Brownian motion also gives insight in how one may go about dealing with greater generality.

\appendix
\renewcommand{\theequation}{A.\arabic{equation}}
\setcounter{equation}{0}
\section*{Appendix: Hypergeometric identity}
An identity for the hypergeometric function that has been used twice in the main body of the text is taken from formula 3.665(2) in \cite{Table}. It  states that,  for any $0<|a|<r$ and $\nu>0$,  as
\begin{equation}
     \int_0^{\pi} \frac{\sin^{d-2} \phi }{(a^2+2a r \cos \phi + r^2)^\nu} \dd\phi=  \frac{1}{r^{2\nu}} B\Big(\frac{d-1}{2}, \frac{1}{2}\Big) \,{_2}F_1 \Big(\nu, \nu-\frac{d}{2}+1; \frac{d}{2}; \frac{a^2}{r^2}\Big) ,
     \label{3.665}
\end{equation}

\section*{Acknowledgements}
All three authors are grateful to two anonymous referees for careful reading of an earlier version of this article and asking challenging, insightful questions, which resulted in a significant  improvements to this article.
TS acknowledges support from a Schlumberger Faculty of the Future award.  SP acknowledges support from the Royal Society as a Newton International Fellow Alumnus (AL191032) and UNAM-DGAPA-PAPIIT grant no. IA103220.

\end{document}